\documentclass[a4paper]{article}
\usepackage[inner=30mm,outer=30mm,textheight=225mm]{geometry}
\usepackage{amsmath}
\usepackage{amssymb}
\usepackage{footnote}
\usepackage{graphicx}
\usepackage{psfrag}
\usepackage{theorem}

\theorembodyfont{\upshape}
\newtheorem{theorem}{Theorem}[subsection]
\newtheorem{corollary}[theorem]{Corollary}
\newtheorem{defn}[theorem]{Definition}
\newtheorem{defns}[theorem]{Definitions}
\newtheorem{example}[theorem]{Example}
\newtheorem{lemma}[theorem]{Lemma}

\newtheorem{algorithm}[theorem]{Algorithm}

\newcommand{\qed}{\hfill\rule[-0.5mm]{1.5mm}{3.0mm}}
\newtheorem{proofthm}{Proof}
\newenvironment{proof}{\begin{proofthm}}{\qed \end{proofthm}}

\renewcommand{\epsilon}{\varepsilon}

\newcommand{\mobius}{M\"{o}bius}

\newcommand{\ppirr}{{$\mathbb{P}^2$-irreducible}}
\newcommand{\ppirrty}{{$\mathbb{P}^2$-irreducibility}}

\newcommand{\regina}{{\em Regina}}

\newcommand{\rps}{\mathbb{R}P^3}

\newcommand{\scircle}{S^1}

\newcommand{\sphere}{S^2}
\newcommand{\sss}{S^3}

\title{Face pairing graphs and 3-manifold enumeration}
\author{Benjamin A.~Burton\footnote{The author would like to acknowledge
    the support of the Australian Research Council.}}
\date{November 20, 2003}

\begin{document}

\maketitle

\abstract{The face pairing graph of a 3-manifold triangulation is a
4-valent graph denoting which tetrahedron faces are identified with
which others.  We present a series of properties that must be satisfied
by the face pairing graph of a closed minimal {\ppirr} triangulation.
In addition we present constraints upon the combinatorial structure of
such a triangulation that can be deduced from its face pairing graph.
These results are then applied to the enumeration of closed minimal
{\ppirr} 3-manifold triangulations, leading to a significant improvement
in the performance of the enumeration algorithm.  Results are offered
for both orientable and non-orientable triangulations.}

\section{Introduction} \label{s-intro}

When studying 3-manifold topology, it is frequently the case that one
has relatively few examples from which to form conjectures, test
hypotheses and draw inspiration.  For this reason a census of all
3-manifold triangulations of a specific type is a useful reference.

Censuses of 3-manifold triangulations have been available in the
literature since 1989 when Hildebrand and Weeks \cite{cuspedcensusold}
published a census of cusped hyperbolic 3-manifold triangulations of up to
five tetrahedra.  This hyperbolic census was later extended to seven
tetrahedra by Callahan, Hildebrand and Weeks \cite{cuspedcensus} and to
closed hyperbolic 3-manifolds by Hodgson and Weeks
\cite{closedhypcensus}.

Matveev \cite{matveev6} in 1998 presented a census of closed orientable
3-manifolds of up to six tetrahedra, extended to seven
tetrahedra by Ovchinnikov (though his results do not appear to be
publicly available) and then to nine tetrahedra by Martelli and
Petronio \cite{italian9}.  Closed non-orientable 3-manifolds were
were tabulated by Amendola and Martelli \cite{italian-nor6} in 2002 for
at most six tetrahedra, extended to seven tetrahedra by Burton
\cite{burton-nor7}.

Each census listed above except for the closed non-orientable census of
Amendola and Martelli \cite{italian-nor6} relies upon a computer
enumeration of 3-manifold triangulations.  Such enumerations typically
involve a brute-force search through possible gluings of tetrahedra and
are thus exceptionally slow, as evidenced by the limits of current
knowledge described above.

As a result the algorithms for enumerating 3-manifold triangulations
generally use mathematical techniques to restrict this brute-force
search and thus improve the overall performance of the enumeration
algorithm.  We present here a series of such techniques based upon face
pairing graphs.  In particular, we derive properties that must be
satisfied by the face pairing graphs of minimal triangulations.
Conversely we also examine structural properties of minimal triangulations
that can be derived from their face pairing graphs.

In the remainder of this section we outline our core assumptions and
the basic ideas behind face pairing graphs.  Section~\ref{s-prelim}
presents combinatorial properties of minimal triangulations that are
required for later results.  We return to face pairing graphs in
Section~\ref{s-main} in which the main theorems are proven.  Finally
Section~\ref{s-enumeration} describes how these results can be
assimilated into the enumeration algorithm and presents some empirical
data.  Appendix~\ref{a-smalltri} illustrates a handful of very
small triangulations that are referred to throughout earlier sections.

\subsection{Assumptions}

As in previous censuses of closed 3-manifold triangulations
\cite{italian-nor6,burton-nor7,italian9,matveev6}, we
restrict our attention to triangulations satisfying the following
properties:
\begin{itemize}
    \item {\em Closed:}  The triangulation is of a closed
    3-manifold.  In particular it has no boundary faces or ideal
    vertices.
    \item {\em \ppirr:}  The underlying 3-manifold has no
    embedded two-sided projective planes, and furthermore every
    embedded 2-sphere bounds a ball.
    \item {\em Minimal:}  The underlying 3-manifold cannot be
    triangulated using strictly fewer tetrahedra.
\end{itemize}

The assumptions of {\ppirrty} and minimality allow us to reduce the
number of triangulations produced in a census to manageable levels.
{\ppirrty} restricts our attention to the simplest 3-manifolds from which
more complex 3-manifolds can be constructed, and minimality restricts
our attention to the simplest triangulations of these 3-manifolds.  It
is worth noting that minimal triangulations are often well structured
(as seen for instance in \cite{burton-nor7}), making them particularly
suitable for studying the 3-manifolds that they represent.

\subsection{Face Pairing Graphs}

A face pairing graph is simply a convenient visual tool for describing
which tetrahedron faces are joined to which others in a 3-manifold
triangulation.  For our purposes we allow graphs to be
{\em multigraphs}, i.e., we allow graphs to contain loops (edges
that join a vertex to itself) and multiple edges (different edges that
join the same pair of vertices).

\begin{defn}[Face Pairing Graph] \label{d-facepairinggraph}
    Let $T$ be a 3-manifold triangulation formed from $n$ tetrahedra,
    and let these tetrahedra be labelled $1,2,\ldots,n$.  The
    {\em face pairing graph} of $T$ is the multigraph on the $n$ vertices
    $1,2,\ldots,n$ constructed as follows.

    Beginning with an empty graph,
    for each pair of tetrahedron faces that are identified in $T$ we
    insert an edge joining the corresponding two graph vertices.
    Specifically, if a face of tetrahedron $i$ is
    identified with a face of tetrahedron $j$ (where $i$ and $j$ may be
    equal) then we insert an edge joining vertices $i$ and $j$.
\end{defn}

\begin{example}
    Consider the two-tetrahedron triangulation of the product space
    $\sphere \times \scircle$ illustrated in Figure~\ref{fig-s2xs1}.
    The tetrahedra are labelled 1 and 2.  The back two faces of
    tetrahedron~1 are identified with each other and the back
    two faces of tetrahedron~2 are
    identified with each other.  The front two faces of
    tetrahedron~1 are each identified with one of the front two
    faces of tetrahedron~2.

    \begin{figure}[htb]
    \psfrag{1}{{\small $1$}} \psfrag{2}{{\small $2$}}
    \centerline{\includegraphics[scale=0.7]{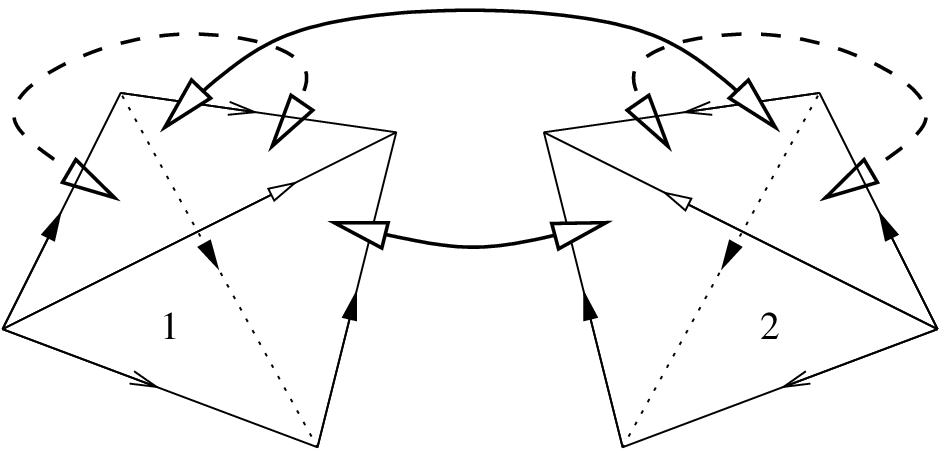}}
    \caption{A two-tetrahedron triangulation of $\sphere \times \scircle$}
    \label{fig-s2xs1}
    \end{figure}

    The corresponding face pairing graph is illustrated in
    Figure~\ref{fig-graphs2xs1}.  Vertex~1 is joined to itself by a
    loop, vertex~2 is joined to itself by a loop and vertices~1 and~2
    are joined together by two distinct edges.

    \begin{figure}[htb]
    \psfrag{1}{{\small $1$}} \psfrag{2}{{\small $2$}}
    \centerline{\includegraphics[scale=0.7]{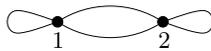}}
    \caption{The face pairing graph corresponding to Figure~\ref{fig-s2xs1}}
    \label{fig-graphs2xs1}
    \end{figure}
\end{example}

The following properties of face pairing graphs are immediate.

\begin{lemma} \label{l-graphprops}
    If $T$ is a closed 3-manifold triangulation formed from $n$ tetrahedra,
    then the face pairing graph of $T$ is a connected multigraph whose
    vertices each have degree 4.
\end{lemma}

\begin{proof}
    Connectedness is immediate since two tetrahedra of $T$ are adjacent
    if and only if the corresponding two graph vertices are adjacent.
    Since triangulation $T$ is closed, each of the four faces of
    every tetrahedron is identified with some other tetrahedron face in
    $T$, and so every graph vertex has degree 4.
\end{proof}

\begin{figure}[htb]
\centerline{\includegraphics[scale=0.7]{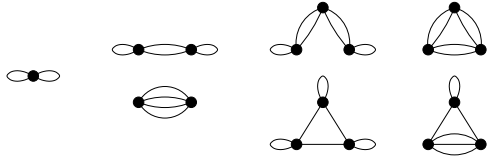}}
\caption{All connected 4-valent multigraphs on $\leq 3$ vertices}
\label{fig-facegraphs123}
\end{figure}

For illustration, the set of all multigraphs on at most three vertices
satisfying the properties of Lemma~\ref{l-graphprops} is presented in
Figure~\ref{fig-facegraphs123}.

\section{Properties of Minimal Triangulations} \label{s-prelim}

In this section we pause to prove a series of basic properties of
minimal triangulations.  We return to face pairing graphs in
Section~\ref{s-main} where the properties described here are used
in proving the main results.

Note that the properties listed in this section are useful in
their own right.  Most of these properties are straightforward to test
by computer, and many can be tested on triangulations that are
only partially constructed.  Thus these properties are easily incorporated into
algorithms for enumerating 3-manifold triangulations.

The properties listed in this section are divided into two
categories.  Section~\ref{s-lowdegedges} presents properties relating to
edges of low degree in a triangulation and
Section~\ref{s-facesubcomplexes} presents properties relating to
subcomplexes formed by faces of a triangulation.

\subsection{Low Degree Edges} \label{s-lowdegedges}

Eliminating edges of low degree whilst enumerating 3-manifold
triangulations is a well-known technique.  The hyperbolic censuses of
Callahan, Hildebrand and Weeks \cite{cuspedcensus,cuspedcensusold} use
results regarding low degree edges in hyperbolic triangulations.
The censuses of Matveev \cite{matveev6} and Martelli and Petronio
\cite{italian9} use similar results for closed orientable
triangulations, which are also proven by Jaco and Rubinstein
\cite{0-efficiency} as consequences of their theory of 0-efficiency.
Here we prove a similar set of results for closed minimal {\ppirr}
triangulations, both orientable and non-orientable.

\begin{lemma}[Degree Three Edges] \label{l-pruneedgedeg3}
    No closed minimal triangulation has an edge of degree three
    that belongs to three distinct tetrahedra.
\end{lemma}

\begin{figure}[htb]
\centerline{\includegraphics[scale=0.7]{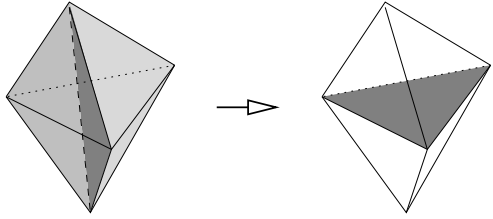}}
\caption{A 3-2 Pachner move}
\label{fig-eltmove32}
\end{figure}

\begin{proof}
    If a triangulation contains an edge of degree three belonging to three
    distinct tetrahedra, a 3-2 Pachner move can be applied to that edge.
    This move involves replacing these three tetrahedra with a pair of
    tetrahedra adjacent along a single face, as illustrated in
    Figure~\ref{fig-eltmove32}.  The resulting triangulation represents
    the same 3-manifold but contains one fewer tetrahedron, and hence the
    original triangulation cannot be minimal.
\end{proof}

The condition of Lemma~\ref{l-pruneedgedeg3} requiring the degree three
edge to belong to three distinct tetrahedra is indeed necessary.  There
do exist closed minimal {\ppirr} triangulations containing edges of degree
three that belong to only two tetrahedra (by representing multiple
edges of one of these tetrahedra).  Examples of such triangulations are
described in detail by Jaco and Rubinstein \cite{0-efficiency}.

Disqualifying edges of degree one or two is more difficult than
disqualifying edges of degree three.  To assist in this task we present
the following simple result.

\begin{lemma} \label{l-embeddedsurfaces}
    No closed minimal {\ppirr} triangulation with $\geq 3$
    tetrahedra contains an embedded non-separating 2-sphere or an
    embedded projective plane.
\end{lemma}

\begin{proof}
    By {\ppirrty}, any embedded 2-sphere
    bounds a ball and is thus separating.  Consider then an embedded
    projective plane $P$.  By {\ppirrty} again we see that
    $P$ is one-sided, and so a regular neighbourhood of $P$ has 2-sphere
    boundary and is in fact $\rps$ with a ball removed.

    Using {\ppirrty} once more it follows that this 2-sphere boundary
    must bound a ball on the side away from $P$
    and so our entire triangulation is of the 3-manifold $\rps$.  We
    observe that $\rps$ can be triangulated using two tetrahedra
    as seen in Appendix~\ref{a-smalltri} and so our triangulation
    cannot be minimal.
\end{proof}

\begin{lemma}[Degree Two Edges] \label{l-pruneedgedeg2}
    No closed minimal {\ppirr} triangulation with $\geq 3$ tetrahedra
    contains an edge of degree two.
\end{lemma}

\begin{proof}
    Let $e$ be an edge of degree two in a closed minimal {\ppirr}
    triangulation with $\geq 3$ tetrahedra.
    If $e$ belongs to only one tetrahedron then it must appear as
    two distinct edges of that tetrahedron.
    Figure~\ref{fig-censusdegree2} lists the three possible arrangements
    in which this is possible.

    \begin{figure}[htb]
    \psfrag{I}{{\small I}} \psfrag{II}{{\small II}} \psfrag{III}{{\small III}}
    \psfrag{e}{{\small $e$}}
    \centerline{\includegraphics[scale=0.7]{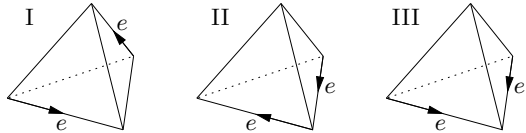}}
    \caption{An edge of degree two belonging to only one tetrahedron}
    \label{fig-censusdegree2}
    \end{figure}

    In case~I, edge $e$ lies in all four faces of the tetrahedron.
    However, since $e$ has degree two it can only belong to two faces of
    the overall triangulation.  Thus these four faces are identified in
    pairs and the triangulation cannot have more than one tetrahedron.

    In cases~II and~III, the bottom face of the tetrahedron contains edge
    $e$ twice and must therefore be identified with some other face of the
    tetrahedron also containing edge $e$ twice.  There is however no
    other such face and so neither of these cases can occur.

    \begin{figure}[htb]
    \psfrag{e}{{\small $e$}} \psfrag{g}{{\small $g$}} \psfrag{h}{{\small $h$}}
    \psfrag{A}{{\small $A$}} \psfrag{B}{{\small $B$}}
    \psfrag{C}{{\small $C$}} \psfrag{D}{{\small $D$}}
    \centerline{\includegraphics[scale=0.7]{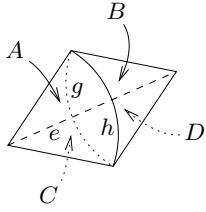}}
    \caption{An edge of degree two belonging to two tetrahedra}
    \label{fig-censusdegree2-2}
    \end{figure}

    Thus $e$ must belong to two distinct tetrahedra as depicted in
    Figure~\ref{fig-censusdegree2-2}.  Consider edges $g$ and $h$ as
    marked in the diagram.  If these edges are identified, the disc
    between them forms either a 2-sphere or a projective plane.  If it
    forms a projective plane then this projective plane is embedded,
    contradicting Lemma~\ref{l-embeddedsurfaces}.

    If the disc between edges $g$ and $h$ forms a 2-sphere, we claim
    that this 2-sphere must be separating.  If the 2-sphere
    does not intersect itself then it is separating by a direct
    application of Lemma~\ref{l-embeddedsurfaces}.  Otherwise the only way in
    which the 2-sphere can intersect itself is if the two vertices at
    the endpoints of $g$ and $h$ are identified.  In this case we deform
    the 2-sphere slightly by pulling it away from these vertices within the
    3-manifold, resulting in an embedded 2-sphere which from
    Lemma~\ref{l-embeddedsurfaces} is again separating.

    So assuming this disc forms a separating 2-sphere, cutting along it
    splits the underlying 3-manifold into a connected sum decomposition.
    Since the triangulation is {\ppirr}, one
    side of this disc must bound a ball; without loss of generality
    let it be the side containing faces $A$ and $C$.  The triangulation
    can then be simplified without changing the underlying 3-manifold by
    removing this ball and flattening the disc to a single edge, converting
    Figure~\ref{fig-censusdegree2-2} to a triangular pillow bounded by faces
    $B$ and $D$ as illustrated in Figure~\ref{fig-censusdegree2-split}.

    \begin{figure}[htb]
    \psfrag{e}{{\small $e$}} \psfrag{g}{{\small $g$}} \psfrag{h}{{\small $h$}}
    \psfrag{gh}{{\small $g,h$}}
    \psfrag{A}{{\small $A$}} \psfrag{B}{{\small $B$}}
    \psfrag{C}{{\small $C$}} \psfrag{D}{{\small $D$}}
    \centerline{\includegraphics[scale=0.7]{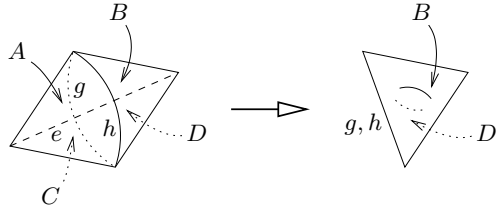}}
    \caption{Simplifying the triangulation by removing a ball}
    \label{fig-censusdegree2-split}
    \end{figure}

    If faces $B$ and $D$ are identified, the 3-manifold
    must be either the 3-sphere or $L(3,1)$ (the two spaces obtainable
    by identifying the faces of a triangular pillow).  Our original
    triangulation is therefore non-minimal since each of these spaces can be
    realised with $\leq 2$ tetrahedra as seen in Appendix~\ref{a-smalltri}.
    If faces $B$ and $D$ are not identified, the entire pillow can be
    flattened to a single face and our original triangulation has been reduced
    by two tetrahedra, once more showing it to be non-minimal.

    The only case remaining is that in which edges $g$ and $h$ are not
    identified.  In this case we may flatten the disc between edges $g$
    and $h$ to a single edge, converting Figure~\ref{fig-censusdegree2-2} to a
    pair of triangular pillows as illustrated in the central diagram of
    Figure~\ref{fig-censusdegree2-pillows}.

    \begin{figure}[htb]
    \psfrag{e}{{\small $e$}} \psfrag{g}{{\small $g$}} \psfrag{h}{{\small $h$}}
    \psfrag{A}{{\small $A$}} \psfrag{B}{{\small $B$}}
    \psfrag{C}{{\small $C$}} \psfrag{D}{{\small $D$}}
    \centerline{\includegraphics[scale=0.7]{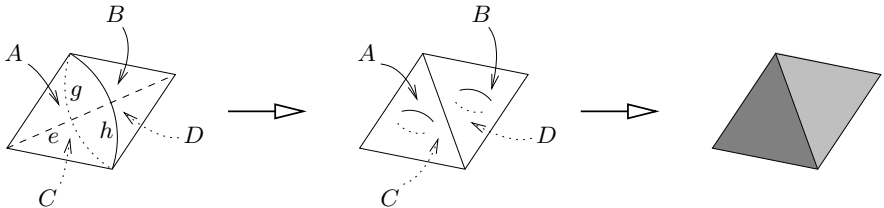}}
    \caption{Flattening the region about an edge of degree two}
    \label{fig-censusdegree2-pillows}
    \end{figure}

    Each of these pillows may then be flattened to a face as illustrated
    in the right hand diagram of Figure~\ref{fig-censusdegree2-pillows}.
    The underlying 3-manifold is only changed if we attempt to flatten a
    pillow whose top and bottom faces are identified, in which case
    our 3-manifold must again be $\sss$ or $L(3,1)$ and so our original
    triangulation is non-minimal.  Otherwise we have reduced our
    triangulation by two tetrahedra, once more a contradiction to minimality.
\end{proof}

\begin{lemma}[Degree One Edges] \label{l-pruneedgedeg1}
    No closed minimal {\ppirr} triangulation with $\geq 3$ tetrahedra
    contains an edge of degree one.
\end{lemma}

\begin{proof}
    The only way of creating an edge of degree one in a 3-manifold
    triangulation is to fold two faces of a tetrahedron together around
    the edge between them, as illustrated in
    Figure~\ref{fig-censusedgedeg1} where $e$ is the edge of degree one.

    \begin{figure}[htb]
    \psfrag{e}{{\small $e$}}
    \centerline{\includegraphics[scale=0.7]{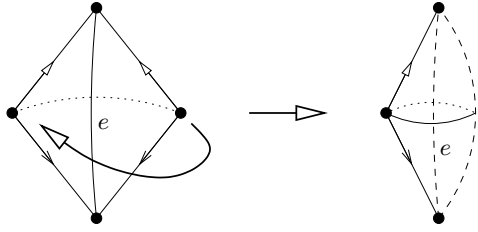}}
    \caption{An edge of degree one}
    \label{fig-censusedgedeg1}
    \end{figure}

    Since our triangulation has $\geq 3$ tetrahedra, the two remaining
    faces of this tetrahedron cannot be identified with each other.
    Let the upper face then be joined to some other tetrahedron as illustrated
    in the left hand diagram of Figure~\ref{fig-census21sum}.  Let $g$ and
    $h$ be the two edges of this new tetrahedron running from $A$ to
    $C$ as marked in the diagram.

    \begin{figure}[htb]
    \psfrag{A}{{\small $A$}} \psfrag{B}{{\small $B$}}
    \psfrag{C}{{\small $C$}} \psfrag{D}{{\small $D$}}
    \psfrag{e}{{\small $e$}} \psfrag{g}{{\small $g$}} \psfrag{h}{{\small $h$}}
    \psfrag{gh}{{\small $g,h$}}
    \psfrag{(M)}{{\small ($M$)}} \psfrag{(N)}{{\small ($N$)}}
    \centerline{\includegraphics[scale=0.7]{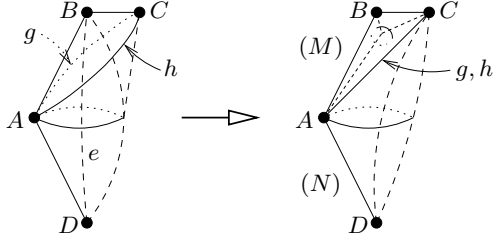}}
    \caption{Splitting apart a connected sum}
    \label{fig-census21sum}
    \end{figure}

    If edges $g$ and $h$ are identified, the disc between them forms
    either a 2-sphere or a projective plane.  Assuming our triangulation
    is minimal and {\ppirr}, the same argument used in the proof of
    Lemma~\ref{l-pruneedgedeg2} shows that this disc must be a
    separating 2-sphere.

    In this case, cutting along the disc between edges $g$ and $h$ therefore
    splits our 3-manifold into two separate pieces each with 2-sphere
    boundary.  Flattening this disc to a single edge effectively fills
    in each of these 2-sphere boundary components with balls, producing
    closed 3-manifolds $M$ and $N$ for which our original 3-manifold is
    the connected sum $M \# N$.

    This procedure is illustrated in the right hand diagram of
    Figure~\ref{fig-census21sum}, where $M$ includes the portion above
    edge $AC$ (including vertex $B$) and $N$ includes the portion below
    edge $AC$ (including vertex $D$).  Note that the portion of the
    diagram between vertices $A$, $C$ and $B$ forms a triangular pillow,
    which can be retriangulated using two tetrahedra with a new internal
    vertex as illustrated.  The portion between vertices $A$, $C$ and
    $D$ can be retriangulated using a single tetrahedron with a new
    internal edge of degree one, much like our original
    Figure~\ref{fig-censusedgedeg1}.

    Since our original 3-manifold is {\ppirr} but is also expressible as
    the connected sum $M \# N$, it follows that one of $M$ or $N$ is a
    3-sphere and the other is a new triangulation of this original
    3-manifold.
    If $M$ is the 3-sphere, we see that $N$ is a triangulation of the
    original 3-manifold formed from one fewer tetrahedron.  If $N$ is the
    3-sphere then $M$ is a new triangulation of our original 3-manifold
    with no change in the number of tetrahedra.  This new triangulation
    however contains three edges of degree two within the triangular
    pillow, and so from Lemma~\ref{l-pruneedgedeg2} it cannot be formed
    from the smallest possible number of tetrahedra.  Either way we see that
    our original triangulation cannot be minimal.

    The only remaining case is that in which edges $g$ and $h$ are not
    identified at all.  Here we may flatten the disc between $g$ and $h$
    to an edge without altering the underlying 3-manifold, as
    illustrated in the central diagram of Figure~\ref{fig-census21}.
    The region between vertices $A$, $C$ and $D$ may be retriangulated
    using a single tetrahedron as before.

    \begin{figure}[htb]
    \psfrag{A}{{\small $A$}} \psfrag{B}{{\small $B$}}
    \psfrag{C}{{\small $C$}} \psfrag{D}{{\small $D$}}
    \psfrag{e}{{\small $e$}} \psfrag{g}{{\small $g$}} \psfrag{h}{{\small $h$}}
    \centerline{\includegraphics[scale=0.7]{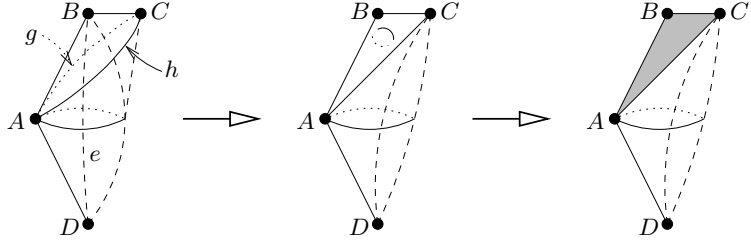}}
    \caption{Flattening a region close to an edge of degree one}
    \label{fig-census21}
    \end{figure}

    In this case we once more observe that the region between
    vertices $A$, $C$ and $B$ becomes a triangular pillow.
    Furthermore, it is impossible to identify the two
    faces bounding this pillow under any rotation or reflection
    without identifying edges $g$ and $h$ as
    a result.  Hence the two faces bounding this pillow are distinct and
    we may flatten the pillow to a face as illustrated in the right-hand
    diagram of Figure~\ref{fig-census21} with no effect upon the
    underlying 3-manifold.  The result is a new triangulation of our
    original 3-manifold with one fewer tetrahedron, and so again our
    original triangulation is non-minimal.
\end{proof}

\subsection{Face Subcomplexes} \label{s-facesubcomplexes}

In this section we identify a variety of simple structures formed from
faces within a closed triangulation whose presence indicates that the
triangulation cannot be both minimal and {\ppirr}.

For orientable triangulations, slightly restricted versions of these
results can be seen as a consequence of Jaco and Rubinstein's theory of
0-efficiency \cite{0-efficiency}.  Here
we prove these results using more elementary means, with generalisations
in the orientable case and extensions to the non-orientable case.

\begin{lemma}[$L(3,1)$ Spines] \label{l-l31spine}
    Let $T$ be a closed minimal {\ppirr} triangulation containing $\geq 3$
    tetrahedra.  Then no face of $T$ has all three of its edges
    identified with each other in the same direction around the face, as
    illustrated in Figure~\ref{fig-censusl31spine}.

    \begin{figure}[htb]
    \centerline{\includegraphics[scale=0.7]{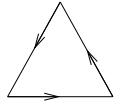}}
    \caption{A face with three identified edges}
    \label{fig-censusl31spine}
    \end{figure}

    In Jaco and Rubinstein's categorisation of face types within a
    triangulation \cite{0-efficiency}, faces of this type are referred
    to as $L(3,1)$ spines for reasons that become apparent in the proof
    below.
\end{lemma}

\begin{proof}
    Consider a regular neighbourhood of the face in question.  This
    regular neighbourhood can be expressed as a triangular prism, where
    the original face slices through the centre of the prism and the
    half-rectangles surrounding this central face are identified in pairs.
    Such structures can be seen in Figure~\ref{fig-censusl31prism}, which
    depicts the two possible neighbourhoods that can be formed in this
    way up to rotation and reflection.  In each diagram the original face
    with the three identified edges is shaded.

    \begin{figure}[htb]
    \psfrag{A}{{\small $A$}} \psfrag{B}{{\small $B$}} \psfrag{C}{{\small $C$}}
    \psfrag{D}{{\small $D$}} \psfrag{E}{{\small $E$}} \psfrag{F}{{\small $F$}}
    \psfrag{X}{{\small $X$}} \psfrag{Y}{{\small $Y$}} \psfrag{Z}{{\small $Z$}}
    \centerline{\includegraphics[scale=0.7]{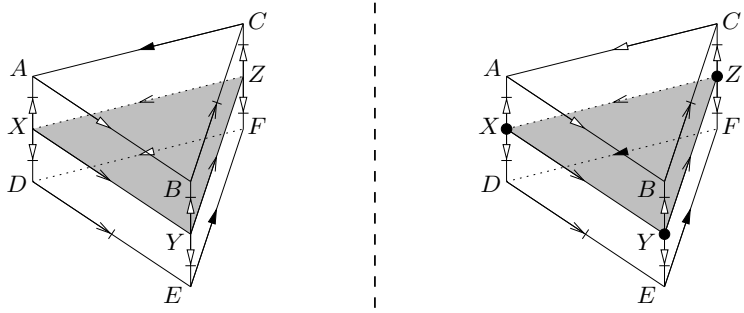}}
    \caption{Possible regular neighbourhoods of a face with three
        identified edges}
    \label{fig-censusl31prism}
    \end{figure}

    In the left hand diagram the neighbourhood is orientable, with
    faces {\em XYBA} and {\em ZXDF} identified, faces {\em YZCB} and
    {\em XYED} identified and faces {\em ZXAC} and {\em YZFE} identified.
    It can be seen that the boundary of this prism is a 2-sphere
    and that the body of the prism forms the lens space $L(3,1)$ with a ball
    removed.  Thus our original triangulation represents the connected
    sum $M \# L(3,1)$ for some $M$, whereby {\ppirrty} implies that it must
    represent $L(3,1)$ itself.  Since $L(3,1)$ can be triangulated using
    only two tetrahedra (see Appendix~\ref{a-smalltri}), it follows that
    triangulation $T$ cannot be minimal.

    The neighbourhood illustrated in the right hand diagram of
    Figure~\ref{fig-censusl31prism} is non-orientable, with
    faces {\em XYBA} and {\em ZXAC} identified, faces {\em YZCB} and
    {\em XYED} identified and faces {\em ZXDF} and {\em YZFE} identified.
    Note that points $X$, $Y$ and $Z$ are all identified as a single vertex
    in the triangulation, marked by a black circle in the diagram.
    A closer examination shows the link of this
    vertex to be a Klein bottle, not a 2-sphere.  Therefore in this case
    triangulation $T$ cannot represent a closed 3-manifold at all.
\end{proof}

Before proceeding further, we recount the following
result regarding vertices within a minimal triangulation.

\begin{lemma} \label{l-multvtx}
    Let $M$ be a closed {\ppirr} 3-manifold that is not $\sss$, $\rps$
    or the lens space $L(3,1)$.  Then every minimal triangulation of $M$
    has precisely one vertex.
\end{lemma}

\begin{proof}
    This result is proven explicitly by Jaco and Rubinstein
    \cite{0-efficiency} for orientable 3-manifolds.
    In the general case, Martelli and Petronio \cite{italian-decomp}
    prove an equivalent result regarding special spines of 3-manifolds.
\end{proof}

The following results eliminate cones within a triangulation formed from
a small number of faces.  These can be seen as generalisations
of the low degree edge results of Section~\ref{s-lowdegedges}.

\begin{lemma}[Two-Face Cones] \label{l-doubleconefaces}
    Let $T$ be a closed {\ppirr} triangulation containing $\geq 3$
    tetrahedra.  Suppose that two distinct faces of $T$ are joined
    together along their edges to form a cone as illustrated in
    Figure~\ref{fig-twofacecone}.  Note that it is irrelevant whether or
    not any of the four edges of this cone (two internal edges and two
    boundary edges) are identified with each other in the overall
    triangulation.

    \begin{figure}[htb]
    \centerline{\includegraphics[scale=0.7]{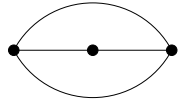}}
    \caption{Two faces joined together to form a cone}
    \label{fig-twofacecone}
    \end{figure}

    Since the two faces of this cone are distinct, if the interior of the cone
    intersects itself then these self-intersections
    can only take place along the two interior edges.
    In addition then, suppose that any such self-intersections are
    tangential and not transverse, as illustrated in
    Figure~\ref{fig-transverse}.
    If all of these conditions are met then $T$ cannot be a minimal
    triangulation.

    \begin{figure}[htb]
    \psfrag{Tangential}{{\small Tangential}}
    \psfrag{Transverse}{{\small Transverse}}
    \centerline{\includegraphics[scale=0.6]{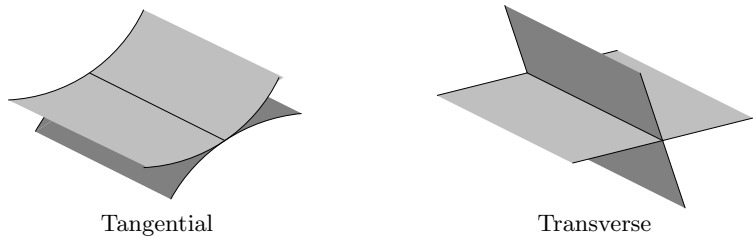}}
    \caption{Tangential and transverse intersections}
    \label{fig-transverse}
    \end{figure}
\end{lemma}

\begin{proof}
    A similar result is proven by Martelli and Petronio in the
    orientable case using special spines of 3-manifolds \cite{italian9}.
    We mirror a key construction of their proof for use in the general
    case here.

    Before beginning, we observe that if triangulation $T$ represents one of
    the 3-manifolds $\sss$, $\rps$ or $L(3,1)$ then $T$ cannot be
    minimal, since each of these 3-manifolds can be triangulated using
    one or two tetrahedra as seen in Appendix~\ref{a-smalltri}.  We
    assume then that $T$ is a triangulation of some other 3-manifold.

    Let the two faces of the cone be $F$ and $G$ with
    interior vertex $v$, as illustrated in the left hand
    diagram of Figure~\ref{fig-twofaceconeexpand}.  Note that this cone forms
    a disc which is by necessity two-sided, and furthermore recall that
    all self-intersections of the interior of the cone are tangential.
    We can therefore expand faces $F$ and
    $G$ into two new tetrahedra as illustrated
    in the right hand diagram of Figure~\ref{fig-twofaceconeexpand}.
    Vertex $v$ is pulled apart into two distinct vertices $v_1$ and $v_2$,
    joined by a new internal edge $e$ of degree two.  Faces $F$ and $G$
    are similarly pulled apart into four distinct
    faces $F_1$, $F_2$, $G_1$ and $G_2$.

    \begin{figure}[htb]
    \psfrag{v}{{\small $v$}}
    \psfrag{v1}{{\small $v_1$}} \psfrag{v2}{{\small $v_2$}}
    \psfrag{F}{{\small $F$}}
    \psfrag{F1}{{\small $F_1$}}
    \psfrag{F2}{{\small $F_2$}}
    \psfrag{G}{{\small $G$}}
    \psfrag{G1}{{\small $G_1$}}
    \psfrag{G2}{{\small $G_2$}}
    \psfrag{e}{{\small $e$}}
    \psfrag{g}{{\small $g$}} \psfrag{h}{{\small $h$}}
    \centerline{\includegraphics[scale=0.6]{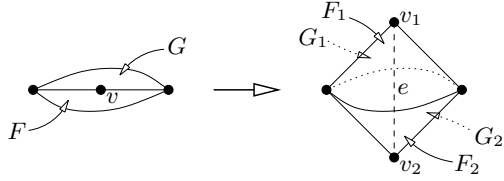}}
    \caption{Expanding a two-face cone into two tetrahedra}
    \label{fig-twofaceconeexpand}
    \end{figure}

    Call the resulting triangulation $T'$.  If the original triangulation $T$
    is formed from $t$ tetrahedra then $T'$ is formed from precisely
    $t+2$ tetrahedra.  Furthermore, $T'$ contains at least two vertices since
    $v_1$ and $v_2$ are distinct.  At this point we attempt to simplify $T'$
    using a similar set of operations to those used in the proof of
    Lemma~\ref{l-pruneedgedeg2}.  These operations are illustrated in
    Figure~\ref{fig-twofaceconesimp}.

    \begin{figure}[htb]
    \psfrag{v}{{\small $v$}}
    \psfrag{v1}{{\small $v_1$}} \psfrag{v2}{{\small $v_2$}}
    \psfrag{F}{{\small $F$}}
    \psfrag{F1}{{\small $F_1$}}
    \psfrag{F2}{{\small $F_2$}}
    \psfrag{G}{{\small $G$}}
    \psfrag{G1}{{\small $G_1$}}
    \psfrag{G2}{{\small $G_2$}}
    \psfrag{e}{{\small $e$}} \psfrag{d}{{\small $\delta$}}
    \psfrag{g}{{\small $g$}} \psfrag{h}{{\small $h$}}
    \centerline{\includegraphics[scale=0.6]{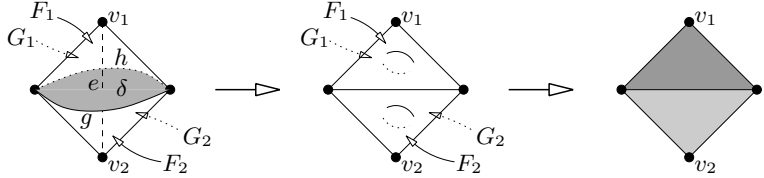}}
    \caption{Simplifying the expanded triangulation $T'$}
    \label{fig-twofaceconesimp}
    \end{figure}

    Let $g$ and $h$ be the edges circling $e$ as marked
    in the left hand diagram of Figure~\ref{fig-twofaceconesimp}, and
    let $\delta$ be the disc that they bound.
    Our first step is to flatten $\delta$ to a single edge as illustrated in
    the central diagram of Figure~\ref{fig-twofaceconesimp}.  If $g$ and
    $h$ are distinct edges of $T'$ then flattening $\delta$ has no effect
    upon the underlying 3-manifold.  Alternatively, if $g$ and $h$ are
    identified in $T'$ then we can use {\ppirrty} as in the
    proofs of Lemmas~\ref{l-pruneedgedeg2} and~\ref{l-pruneedgedeg1}
    to show that the disc $\delta$ must form a separating 2-sphere.  In this
    case, flattening $\delta$ to an edge corresponds to slicing along
    this 2-sphere and filling the resulting boundary components with balls.
    That is, the underlying 3-manifold is broken up into a connected sum
    decomposition.

    At this stage the two new tetrahedra have become two triangular
    pillows.  Each of these pillows can then be flattened to a face as
    illustrated in the right hand diagram of
    Figure~\ref{fig-twofaceconesimp}.  Since all four faces
    $F_1$, $F_2$, $G_1$ and $G_2$ are distinct, this flattening of
    pillows has no effect upon the underlying 3-manifold.

    The final triangulation then has $t$ tetrahedra but more importantly
    has at least two vertices.
    If the underlying 3-manifold was not changed, it
    follows from Lemma~\ref{l-multvtx} that this final triangulation
    is not minimal and so neither is $T$.

    If the underlying 3-manifold was changed however, it was broken into a
    connected sum decomposition.  By {\ppirrty} it follows that one of
    the summands must be a 3-sphere and the other must be the
    original 3-manifold.  In this case we simply throw away the
    3-sphere and obtain a triangulation of the original 3-manifold
    formed from strictly fewer than $t$ tetrahedra, once more showing the
    original triangulation $T$ to be non-minimal.
\end{proof}

\begin{lemma}[One-Face Cones] \label{l-conefaces}
    Let $T$ be a closed minimal {\ppirr} triangulation containing $\geq 3$
    tetrahedra.  Then no single face of $T$ has two of its edges identified
    to form a cone as illustrated in Figure~\ref{fig-facecone}.
    Note that it is irrelevant whether or not the third edge of the face
    is identified with the other two.

    \begin{figure}[htb]
    \psfrag{V}{{\small $v$}}
    \centerline{\includegraphics[scale=0.7]{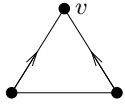}}
    \caption{A face with two edges identified to form a cone}
    \label{fig-facecone}
    \end{figure}
\end{lemma}

\begin{proof}
    Suppose that some face $F$ has two of its edges identified to form a
    cone, where $v$ is the vertex joining both edges as marked in
    Figure~\ref{fig-facecone}.
    Consider the link of vertex $v$ within the 3-manifold triangulation.
    This link must be a 2-sphere formed from triangles, which we denote by $S$.
    Furthermore, face $F$ corresponds to a loop within this vertex link
    as illustrated in the left hand diagram of Figure~\ref{fig-conelink}.

    \begin{figure}[htb]
    \psfrag{I}{{\small I}}
    \psfrag{II}{{\small II}}
    \psfrag{III}{{\small III}}
    \centerline{\includegraphics[scale=0.7]{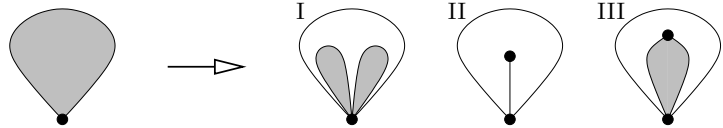}}
    \caption{The inside of a loop in a vertex link}
    \label{fig-conelink}
    \end{figure}

    Since $S$ is a 2-sphere, this loop bounds a disc in the vertex link.
    The inside of
    this disc must appear as one of the three configurations shown on
    the right hand side of Figure~\ref{fig-conelink}, where the
    structures of the shaded regions remain unknown.  Configuration~I
    contains one or more inner loops; we can thus
    apply this same argument to the inner loops and continue deeper
    inside the disc until we eventually arrive at either
    configuration~II or~III.

    In configuration~II we see that the link $S$
    contains a vertex of degree one.  This corresponds to an edge
    of degree one in the 3-manifold triangulation, which from
    Lemma~\ref{l-pruneedgedeg1} cannot occur.

    In configuration~III, the two inner edges of the vertex link
    correspond to two faces $G$ and $H$ of the 3-manifold triangulation
    whose edges
    are identified to form a cone, as illustrated in the left hand
    diagram of Figure~\ref{fig-conebody}.  Furthermore, the triangle
    wrapping around these edges in the vertex link corresponds to a
    tetrahedron wrapping around faces $G$ and $H$
    as seen in the right hand diagram of Figure~\ref{fig-conebody}.

    \begin{figure}[htb]
    \psfrag{G}{{\small $G$}}
    \psfrag{H}{{\small $H$}}
    \psfrag{v}{{\small $v$}}
    \centerline{\includegraphics[scale=0.7]{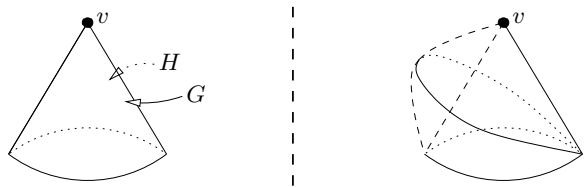}}
    \caption{Faces corresponding to the two inner edges of the vertex link}
    \label{fig-conebody}
    \end{figure}

    If faces $G$ and $H$ are distinct, we call upon
    Lemma~\ref{l-doubleconefaces} to secure the result.
    The tetrahedron wrapped around these
    faces ensures that the required tangentiality condition is met, and
    so Lemma~\ref{l-doubleconefaces} can be
    applied to show that triangulation $T$ is non-minimal.

    \begin{figure}[htb]
    \centerline{\includegraphics[scale=0.7]{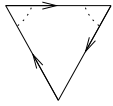}}
    \caption{A face with two corners identified to form a cone}
    \label{fig-conefacel31}
    \end{figure}

    Otherwise $G$ and $H$ must in fact be the same face of the
    3-manifold triangulation, where the two inner edges of the vertex link
    correspond to two different corners of the face.  For two different
    corners to meet in this way, this face must have
    all three edges identified as illustrated in Figure~\ref{fig-conefacel31}.
    Lemma~\ref{l-l31spine} however shows that such a face cannot occur.
\end{proof}

Our next results concern pairs of faces that, when sliced open, combine to
form small spherical boundary components of the resulting triangulation.

\begin{lemma}[Spherical Subcomplexes] \label{l-splittripillow}
    Let $T$ be a closed triangulation with $\geq 3$ tetrahedra.
    Consider two faces $F_1$ and $F_2$ of $T$ that are joined along at
    least one edge.

    Slicing $T$ along $F_1$ and $F_2$ will produce a new (possibly
    disconnected) triangulation with four boundary faces; call this $T'$.
    Note that $T'$ might have vertices whose links are neither spheres
    nor discs, and so might not actually represent one or more 3-manifolds.

    If $T'$ contains multiple boundary components (as opposed to a
    single four-face boundary component) and if one of these boundary
    components is a sphere as illustrated in
    Figure~\ref{fig-censustripillow}, then $T$ cannot be both minimal
    and {\ppirr}.

    \begin{figure}[htb]
    \centerline{\includegraphics[scale=0.7]{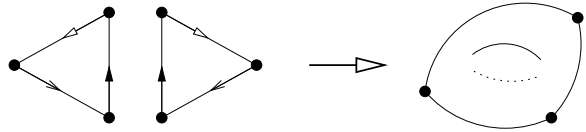}}
    \caption{A sphere formed by identifying the boundaries of two triangles}
    \label{fig-censustripillow}
    \end{figure}

    Note that it does not matter if some vertices of this boundary
    sphere are identified, i.e., if the sphere is pinched at two or
    more points -- this is indeed expected if $T'$ has non-standard
    vertex links as suggested above.  Only the edge identifications
    illustrated in Figure~\ref{fig-censustripillow} are required for the
    conditions of this lemma to be met.
\end{lemma}

\begin{proof}
    Let $T$ be a minimal {\ppirr} triangulation satisfying all of the
    given conditions and consider the triangulation $T'$ as described in
    the lemma statement.  Recall that any boundary component must contain
    an even number of faces.  Since $T'$ has four boundary faces and
    multiple boundary components, it must therefore have precisely
    two boundary components each with
    two faces.  Furthermore, since faces $F_1$ and $F_2$ are adjacent in
    $T$, each of these boundary components must be formed from a copy of
    both faces $F_1$ and $F_2$.

    Let these two boundary components be $\partial_1$ and
    $\partial_2$ where $\partial_1$ is a sphere as illustrated in
    Figure~\ref{fig-censustripillow}.  Let $R_1$ and $R_2$ represent the
    regions of $T'$ just inside boundary components $\partial_1$ and
    $\partial_2$ respectively.  Finally let $S$ be an embedded sphere located
    in $R_1$ just behind the spherical boundary component $\partial_1$.
    This scenario is illustrated in Figure~\ref{fig-censustripillowsetup}.

    \begin{figure}[htb]
    \psfrag{R1}{{\small $R_1$}} \psfrag{R2}{{\small $R_2$}}
    \psfrag{d1}{{\small $\partial_1$}} \psfrag{d2}{{\small $\partial_2$}}
    \psfrag{S}{{\small $S$}}
    \centerline{\includegraphics[scale=0.7]{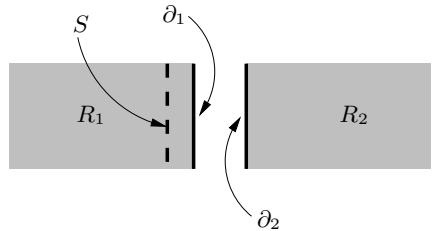}}
    \caption{Regions within the split triangulation $T'$}
    \label{fig-censustripillowsetup}
    \end{figure}

    As noted in the lemma statement, vertices of the boundary sphere
    $\partial_1$ might be identified and thus it may in fact be
    impossible to place the sphere $S$ entirely within region $R_1$.  If
    this is the case then we simply push $S$ outside region $R_1$ in a
    small neighbourhood of each offending vertex.  Examples of the
    resulting sphere $S$ are illustrated for two different cases in
    Figure~\ref{fig-censustripillowpushs}.

    \begin{figure}[htb]
    \psfrag{R1}{{\small $R_1$}}
    \psfrag{d1}{{\small $\partial_1$}}
    \psfrag{S}{{\small $S$}}
    \centerline{\includegraphics[scale=0.7]{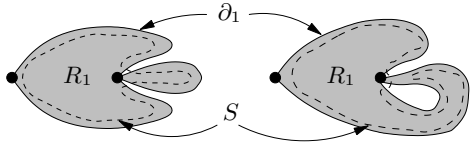}}
    \caption{Pushing sphere $S$ outside region $R_1$ in the neighbourhood
        of a vertex}
    \label{fig-censustripillowpushs}
    \end{figure}

    We can now view $S$ as an embedded 2-sphere in the original
    triangulation $T$.  Note that whenever $S$ is
    pushed outside region $R_1$ as described above, it is simply
    pushed into a small neighbourhood of a vertex in region
    $R_2$.  In particular, since the link of each vertex of $T$ is a
    2-sphere, such operations do not introduce any self-intersections in $S$.

    Since $T$ is {\ppirr} it follows that the embedded 2-sphere
    $S$ must bound a ball in $T$.  We take cases according to whether
    this ball lies on the side of $S$ including region $R_1$ or region $R_2$.

    \begin{itemize}
        \item Suppose that $S$ bounds a ball on the side containing region
        $R_2$.  In this case we can remove the component of $T'$ containing
        region $R_2$ and replace it with a two-tetrahedron
        triangular pillow as
        illustrated in Figure~\ref{fig-censustripillowplug}.  This
        pillow can then be joined to region $R_1$ along the spherical
        boundary $\partial_1$ resulting in a new triangulation $T''$
        whose underlying 3-manifold is identical to that of $T$.

        \begin{figure}[htb]
        \centerline{\includegraphics[scale=0.7]{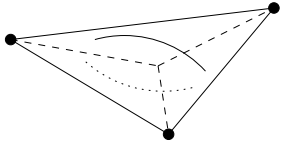}}
        \caption{A replacement two-tetrahedron triangular pillow}
        \label{fig-censustripillowplug}
        \end{figure}

        \item Otherwise $S$ bounds a ball on the side containing region $R_1$.
        In this case we remove the component of $T'$ containing region
        $R_1$ and again replace it with the two-tetrahedron triangular
        pillow illustrated in Figure~\ref{fig-censustripillowplug}.  The
        pillow has a spherical boundary identical to the old boundary
        $\partial_1$ and so it can be successfully dropped in as a
        replacement component without changing the underlying 3-manifold.
        Once more we denote this new triangulation $T''$, where the
        underlying 3-manifold of $T''$ is identical to that of $T$.

        Note that even if the old component of $T'$ containing region $R_1$
        includes a vertex whose link is a multiply-punctured sphere, the
        underlying 3-manifolds of $T$ and $T''$ remain identical.
        Such a scenario is illustrated in the left hand diagram of
        Figure~\ref{fig-censustripillowpunc}, with the replacement
        triangular pillow illustrated in the right hand diagram of the
        same figure.  Since $S$ is pushed away from the vertex
        in question, it can be seen that in both triangulations $T$ and
        $T''$ the sphere $S$ bounds a ball on the side containing region
        $R_1$.  Meanwhile the region on the other side of $S$ remains
        unchanged and so the underlying 3-manifold is preserved.

        \begin{figure}[htb]
        \psfrag{R1}{{\small $R_1$}}
        \psfrag{d1}{{\small $\partial_1$}}
        \psfrag{S}{{\small $S$}}
        \psfrag{T/T'}{{\small $T$, $T'$}} \psfrag{T''}{{\small $T''$}}
        \centerline{\includegraphics[scale=0.7]{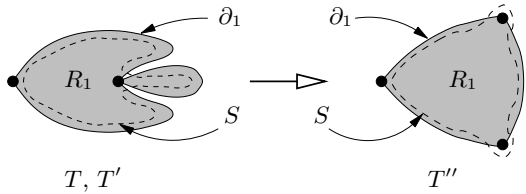}}
        \caption{A case involving a vertex whose link is a
            multiply-punctured sphere}
        \label{fig-censustripillowpunc}
        \end{figure}
    \end{itemize}

    We see then that
    in each of the above cases we have removed at least one tetrahedron
    and inserted two, resulting in a net increase of at most one
    tetrahedron.  However, since the two original faces $F_1$ and $F_2$
    were distinct, we can flatten the triangular pillow of
    Figure~\ref{fig-censustripillowplug} to a single face in the new
    triangulation $T''$.  The resulting
    triangulation has the same underlying 3-manifold as before
    but now uses at least one fewer tetrahedron than the original
    triangulation $T$.  Thus triangulation $T$ cannot be minimal.
\end{proof}

\begin{corollary} \label{c-purgetripillow}
    Let $T$ be a closed triangulation with $\geq 3$ tetrahedra.
    If $T$ contains two faces whose edges are identified to form a
    sphere as illustrated in Figure~\ref{fig-censustripillow} and if the
    three edges of this sphere are distinct (i.e., none are identified
    with each other in the triangulation), then $T$ cannot be both
    minimal and {\ppirr}.
\end{corollary}

\begin{proof}
    Since the three edges of this sphere are distinct, the sphere
    has no self-intersections except possibly at the vertices
    of its two faces.  Combined with the fact that all spheres within
    a 3-manifold are two-sided, we see that slicing $T$ along the two
    faces in question produces two spherical boundary components that are
    each triangulated as illustrated in
    Figure~\ref{fig-censustripillow}.  Thus
    Lemma~\ref{l-splittripillow} can be invoked and $T$ cannot be both
    minimal and {\ppirr}.
\end{proof}

The final result of this section, though not strictly falling into the
category of face subcomplexes, nevertheless helps lighten some of
the case analyses that take place in Section~\ref{s-main}.

\begin{lemma} \label{l-tetringor}
    Let $T$ be a closed minimal {\ppirr} triangulation in which $e$ is
    an edge.  Consider a ring of adjacent tetrahedra
    $\Delta_1,\ldots,\Delta_k$ about edge $e$ as illustrated in the left
    hand diagram of
    Figure~\ref{fig-tetring}.  Note that the same tetrahedron may appear
    multiple times in this list, since edge $e$ may correspond to
    multiple edges of the same tetrahedron.  Let $F$ and $G$ be the two
    faces at either end of this ring, so that faces $F$ and $G$ belong to
    tetrahedra $\Delta_1$ and $\Delta_k$ respectively and are adjacent
    along edge $e$.  These two faces are shaded in the diagram.

    \begin{figure}[htb]
    \psfrag{e}{{\small $e$}}
    \psfrag{d1}{{\small $\Delta_1$}} \psfrag{d2}{{\small $\Delta_2$}}
    \psfrag{dd}{{\small $\Delta_k$}}
    \psfrag{F}{{\small $F$}} \psfrag{G}{{\small $G$}}
    \psfrag{P}{{\small $P$}} \psfrag{Q}{{\small $Q$}}
    \psfrag{R}{{\small $R$}} \psfrag{S}{{\small $S$}}
    \centerline{\includegraphics[scale=0.7]{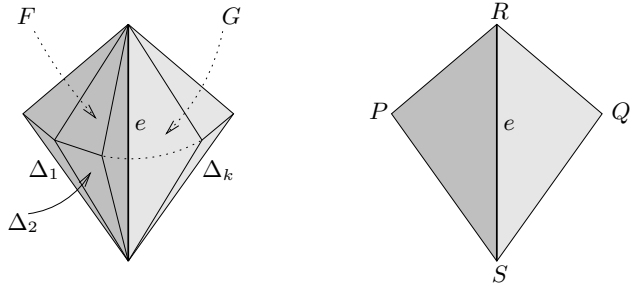}}
    \caption{A ring of tetrahedra about edge $e$}
    \label{fig-tetring}
    \end{figure}

    If faces $F$ and $G$ are identified in triangulation $T$, then this
    identification must be orien\-ta\-tion-pre\-serving within the ring of
    tetrahedra.  That is,
    either the two faces are folded together over edge $e$ or the two
    faces are folded together with a $120^\circ$ twist.

    More specifically, if the vertices of faces $F$ and $G$ are labelled
    $P$, $Q$, $R$ and $S$ as indicated in the right hand diagram of
    Figure~\ref{fig-tetring}, then face {\em PRS} may be identified
    with either {\em QRS}, {\em RSQ} or {\em SQR} (where the order of
    vertices indicates the specific rotation or reflection used).  None
    of the orientation-reversing identifications in which {\em PRS} is
    identified with {\em QSR}, {\em SRQ} or {\em RQS} may be used.
\end{lemma}

\begin{proof}
    If face {\em PRS} is identified with {\em QSR} then edge $e$ is
    identified with itself in reverse and so $T$ cannot be a 3-manifold
    triangulation.  If face {\em PRS} is identified with {\em SRQ} then
    edges {\em PR} and {\em SR} of face $F$ are identified to form a cone,
    contradicting Lemma~\ref{l-conefaces}.  Identifying faces {\em PRS} and
    {\em RQS} similarly produces a cone in contradiction to
    Lemma~\ref{l-conefaces}.
\end{proof}

\section{Main Results} \label{s-main}

Equipped with the preliminary results of Section~\ref{s-prelim}, we can
begin to prove our main results involving face pairing graphs and
minimal triangulations.  These results can be split into two categories.
In Section~\ref{s-mainstruct} we use face pairing graphs to derive
structural properties of the corresponding triangulations.
Section~\ref{s-mainbadgraphs} on the other hand establishes properties
of the face pairing graphs themselves.

\subsection{Structural Properties of Triangulations} \label{s-mainstruct}

Recall from Definition~\ref{d-facepairinggraph} that each edge within a
face pairing graph corresponds to a pair of tetrahedron faces that are
identified in a triangulation.  A triangulation cannot be entirely
reconstructed from its face pairing graph however, since each pair of
faces can be identified in one of six different ways (three possible
rotations and three possible reflections).

In this section we identify certain subgraphs within face pairing
graphs and deduce information about the specific rotations and
reflections under which the corresponding tetrahedron faces are
identified.

A key structure that appears frequently within minimal triangulations
and is used throughout the following results is the layered solid torus.
Layered solid tori have been discussed in a variety of informal contexts
by Jaco and Rubinstein.  They appear in \cite{0-efficiency} and are
treated thoroughly by these authors in \cite{layeredlensspaces}.
Analogous constructs involving special spines of
3-manifolds have been described in detail by Matveev \cite{matveev6} and
by Martelli and Petronio \cite{italianfamilies}.

In order to describe the construction of a layered solid torus we introduce
the process of layering.
Layering is a transformation that, when applied to a triangulation with
boundary, does not change the underlying 3-manifold but does change the
curves formed by the boundary edges of the triangulation.

\begin{defn}[Layering] \label{d-layer}
    Let $T$ be a triangulation with boundary and let $e$ be one of its
    boundary edges.  To {\em layer a tetrahedron on edge $e$}, or just
    to {\em layer on edge $e$}, is to take
    a new tetrahedron $\Delta$, choose two of its faces and identify them
    with the two boundary faces on either side of
    $e$ without twists.  This procedure is
    illustrated in Figure~\ref{fig-layering}.

    \begin{figure}[htb]
    \psfrag{e}{{\small $e$}}
    \psfrag{f}{{\small $f$}}
    \psfrag{T}{{\small $\Delta$}}
    \psfrag{M}{{\small $T$}}
    \centerline{\includegraphics[width=5.5cm]{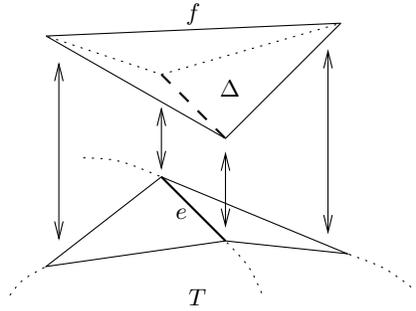}}
    \caption{Layering a tetrahedron on a boundary edge}
    \label{fig-layering}
    \end{figure}
\end{defn}

Note that layering on a boundary edge does not change the underlying
3-manifold; the only effect is to thicken the boundary around edge
$e$.  Note furthermore that edge $e$ is no longer a boundary edge, and
instead edge $f$ (which in general represents a different curve on
the boundary of the 3-manifold) has been added as a new boundary edge.

\begin{defn}[Layered Solid Torus] \label{d-lst}
    A {\em standard layered solid torus} is a triangulation of a solid torus
    formed as follows.  We begin with the {\mobius} band illustrated in
    Figure~\ref{fig-layermobius}, where the two edges marked $e$ are
    identified according to the arrows and where edge $g$ is a boundary edge.
    If we thicken this \mobius\ band slightly, we can imagine it as a
    solid torus with two boundary faces, one on each side of the
    {\mobius} band.

    \begin{figure}[htb]
    \psfrag{e}{{\small $e$}} \psfrag{f}{{\small $f$}} \psfrag{g}{{\small $g$}}
    \centerline{\includegraphics[scale=0.7]{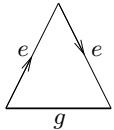}}
    \caption{A one-triangle \mobius\ band}
    \label{fig-layermobius}
    \end{figure}

    We make an initial layering upon one of the boundary edges
    corresponding to edge $e$.  This results in a one-tetrahedron
    triangulation of a solid torus as illustrated in
    Figure~\ref{fig-layermobiusstandard}.
    The new tetrahedron in the illustration is sliced in half; this tetrahedron
    in fact sits on top of the {\mobius} band, runs off the diagram to the
    right and returns from the left to simultaneously sit beneath the
    {\mobius} band.

    \begin{figure}[htb]
    \psfrag{e}{{\small $e$}} \psfrag{g}{{\small $g$}}
    \centerline{\includegraphics[scale=0.7]{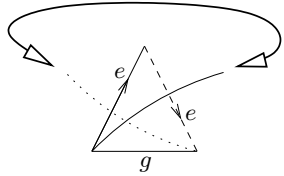}}
    \caption{A one-tetrahedron standard layered solid torus}
    \label{fig-layermobiusstandard}
    \end{figure}

    We then perform some number of additional layerings upon boundary
    edges, one at a time.  We may layer as many times we like or we may
    make no additional layerings at all.  There are thus infinitely many
    different standard layered solid tori that can be constructed.

    Note that a {\em non-standard layered solid torus} can be formed by
    making the initial layering upon edge $g$ of the {\mobius} band
    instead of edge $e$, and that the {\mobius} band itself can be considered
    a {\em degenerate layered solid torus} involving no layerings at all.
    Such structures however are not considered here.
\end{defn}

We can observe that each layered solid torus has two boundary faces
and represents the same underlying 3-manifold, i.e., the solid torus.
What distinguishes the different layered solid tori
is the different patterns of curves that their boundary edges
make upon the boundary torus.

We return now to structures within face pairing graphs.  In particular
we take an interest in chains within face pairing graphs, as defined
below.

\begin{defns}[Chain]
    A {\em chain of length $k$} is the multigraph formed as
    follows.  Take $k+1$ vertices labelled $0,1,2,\ldots,k$ and join
    vertices $i$ and $i+1$ with a double edge for all $0 \leq i \leq k-1$.
    Each of these edges is called an {\em interior edge} of the chain.

    If a loop is added joining vertex 0 to itself the chain becomes a
    {\em one-ended chain}.  If another loop is added joining vertex $k$
    to itself the chain becomes a {\em double-ended chain} (and is now a
    4-valent multigraph).  These loops are called {\em end edges} of
    the chain.
\end{defns}

\begin{example}
    A one-ended chain of length 4 is illustrated in
    Figure~\ref{fig-censuschainoneend}, and a double-ended chain of
    length 3 is illustrated in
    Figure~\ref{fig-censuschain3}.

    \begin{figure}[htb]
    \centerline{\includegraphics[scale=0.7]{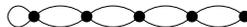}}
    \caption{A one-ended chain in a face pairing graph}
    \label{fig-censuschainoneend}
    \end{figure}

    \begin{figure}[htb]
    \centerline{\includegraphics[scale=0.7]{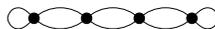}}
    \caption{A double-ended chain in a face pairing graph}
    \label{fig-censuschain3}
    \end{figure}
\end{example}

It can be observed that a one-ended chain of length $k$ is in fact the
face pairing graph of a standard layered solid torus containing $k+1$
tetrahedra.
Our first major result regarding face pairing graphs is a strengthening
of this relationship as follows.

\begin{theorem} \label{t-censuschainoneend}
    Let $T$ be a closed minimal {\ppirr} triangulation with
    $\geq 3$ tetrahedra and face pairing graph $G$.
    If $G$ contains a one-ended chain then the
    tetrahedra of $T$ corresponding to the vertices of this one-ended
    chain form a standard layered solid torus in $T$.
\end{theorem}

\begin{proof}
    We prove this by induction on the chain length.  A one-ended chain of
    length 0 consists of a single vertex with a single end edge, representing
    a single tetrahedron with two of its faces identified.  Using
    Lemma~\ref{l-tetringor} with a ring of just one tetrahedron we see
    that these faces must be identified in an orientation-preserving fashion.

    If these faces
    are simply snapped shut as illustrated in the left hand diagram of
    Figure~\ref{fig-onetetjoin} (faces {\em ABC} and {\em DBC} being
    identified), the edge between them will have
    degree one in the final triangulation which cannot happen
    according to Lemma~\ref{l-pruneedgedeg1}.  Thus
    these faces are identified with a twist as illustrated in the right
    hand diagram of Figure~\ref{fig-onetetjoin} (faces {\em ABC} and
    {\em BCD} being identified), producing a one-tetrahedron standard
    layered solid torus.

    \begin{figure}[htb]
    \psfrag{A}{{\small $A$}} \psfrag{B}{{\small $B$}}
    \psfrag{C}{{\small $C$}} \psfrag{D}{{\small $D$}}
    \centerline{\includegraphics[scale=0.6]{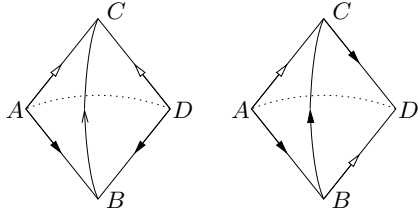}}
    \caption{Different ways of identifying two faces of a single tetrahedron}
    \label{fig-onetetjoin}
    \end{figure}

    Assume then that any one-ended chain of length $k$ must correspond
    to a standard layered solid torus (which will have $k+1$ tetrahedra), and
    consider a one-ended chain of length $k+1$.  This one-ended chain is
    simply a one-ended chain of length $k$ with an extra double edge
    attached to the end, and so by our inductive hypothesis the
    corresponding $k+2$ tetrahedra must form a $(k+1)$-tetrahedron layered
    solid torus with an addition tetrahedron joined to its boundary along
    two faces.

    By symmetry of the two-triangle torus which forms the layered solid
    torus boundary, we can picture the situation as shown in the left
    hand diagram of Figure~\ref{fig-censuslstnextoptions}.
    In this diagram face {\em WXY} of the new tetrahedron is to be
    joined directly to face {\em ABC} of the layered solid torus with
    no rotations or reflections, and face {\em YZW}
    of the new tetrahedron is to be joined to face {\em CDA} but possibly
    with the vertices identified in some different order, i.e., with the
    faces being rotated or reflected before they are identified.
    Applying Lemma~\ref{l-tetringor} to the ring of tetrahedra about
    edge {\em AC}, we see that this identification of faces
    {\em YZW} and {\em CDA} must be orientation-preserving.

    \begin{figure}[htb]
    \centerline{\includegraphics[scale=0.7]{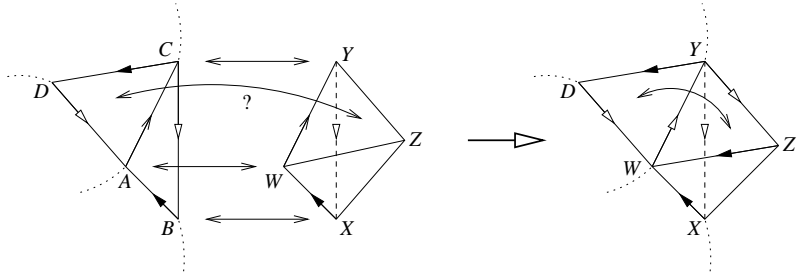}}
    \caption{Attaching a new tetrahedron to a layered solid torus}
    \label{fig-censuslstnextoptions}
    \end{figure}

    If face {\em YZW} is simply folded over onto face {\em CDA}, i.e., the
    vertices are identified in this precise order with no rotation or
    reflection, then we have merely layered the new tetrahedron onto
    edge {\em AC} and so we obtain a larger standard layered solid torus
    with $k+2$ tetrahedra as required.

    Otherwise a rotation must take place, and we may by symmetry
    assume that face {\em YZW} is identified with face {\em ACD}
    as illustrated
    in the right hand diagram of Figure~\ref{fig-censuslstnextoptions}.
    In this diagram we see that each of the two new boundary faces
    ({\em XWZ} and {\em XYZ}) has two edges
    identified to form a cone, in contradiction to Lemma~\ref{l-conefaces}.
\end{proof}

The following result corresponds to double edges in a face pairing
graph, and is particularly helpful in easing the massive case analyses
required for some of the theorems of Section~\ref{s-mainbadgraphs}.

\begin{lemma} \label{l-censusdoubleedge}
    Let $T$ be a closed minimal {\ppirr} triangulation with
    $\geq 3$ tetrahedra.  If two distinct tetrahedra of $T$ are joined
    to each other along two distinct faces then these two face
    identifications must be as illustrated in one of the diagrams of
    Figure~\ref{fig-censusdoubleedgegoodids}.

    \begin{figure}[htb]
    \psfrag{A}{{\small $A$}} \psfrag{B}{{\small $B$}}
    \psfrag{C}{{\small $C$}} \psfrag{D}{{\small $D$}}
    \psfrag{a}{{\small $A'$}} \psfrag{b}{{\small $B'$}}
    \psfrag{c}{{\small $C'$}} \psfrag{d}{{\small $D'$}}
    \centerline{\includegraphics[scale=0.7]{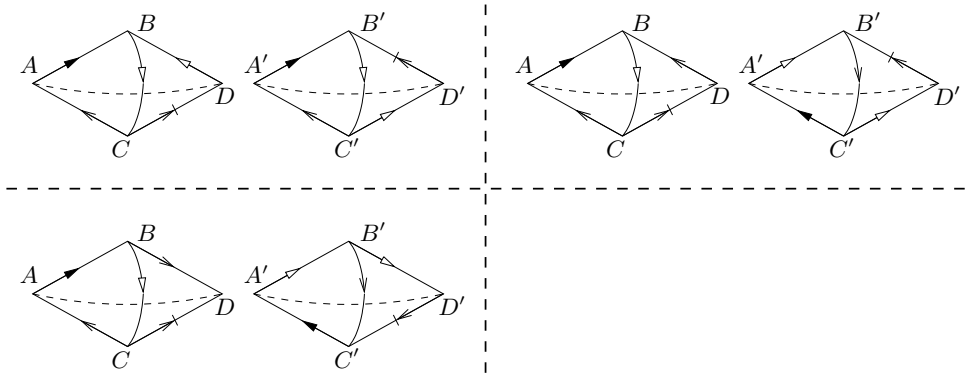}}
    \caption{Allowable ways of joining two tetrahedra along two faces}
    \label{fig-censusdoubleedgegoodids}
    \end{figure}

    Specifically, let the tetrahedra be {\em ABCD} and $A'B'C'D'$ with
    faces {\em ABC} and $A'B'C'$ identified and with faces {\em BCD} and
    $B'C'D'$ identified (though not necessarily with these precise
    vertex identifications; the faces may be rotated or reflected before
    they are identified).
    Then, allowing for the two tetrahedra to be relabelled and/or
    swapped, we must have one of the following three cases.
    \begin{itemize}
        \item Face {\em ABC} is identified with $A'B'C'$ (no rotations
        or reflections take place) and face {\em BCD} is identified with
        $C'D'B'$ (the faces are rotated before being identified);
        this is illustrated in the top left diagram of
        Figure~\ref{fig-censusdoubleedgegoodids}.

        \item Face {\em ABC} is identified with
        $C'A'B'$ and face {\em BCD} is identified with $C'D'B'$
        (both identifications involve a rotation in the same direction);
        this is illustrated in the top right diagram of
        Figure~\ref{fig-censusdoubleedgegoodids}.

        \item Face {\em ABC} is identified with
        $C'A'B'$ (the faces are rotated before being identified)
        and face {\em BCD} is identified with $B'D'C'$ (the faces
        are reflected about a non-horizontal axis); this is illustrated
        in the bottom left diagram of
        Figure~\ref{fig-censusdoubleedgegoodids}.
    \end{itemize}

    Note that the first two cases use consistent orientations for the
    two face identifications, whereas the third case results in a
    non-orientable structure.
\end{lemma}

\begin{proof}
    Consider initially the case in which face {\em ABC} is identified with
    face $A'B'C'$ without reflection or rotation.  Applying
    Lemma~\ref{l-tetringor} to edge {\em AB} we see that the remaining face
    identification must be orientation-preserving.  Thus either face
    {\em BCD} is rotated before being identified with $B'C'D'$,
    resulting in the first case listed in the lemma statement, or
    face {\em BCD} is identified directly with $B'C'D'$, resulting in
    edge {\em BC} having degree two in contradiction to
    Lemma~\ref{l-pruneedgedeg2}.

    Allowing for the tetrahedra to be relabelled and/or swapped, the
    only remaining methods of identifying the two pairs of faces are the
    second and third cases listed in the lemma statement plus the
    additional possibility illustrated in
    Figure~\ref{fig-doubleedgebad}.

    \begin{figure}[htb]
    \psfrag{A}{{\small $A$}} \psfrag{B}{{\small $B$}}
    \psfrag{C}{{\small $C$}} \psfrag{D}{{\small $D$}}
    \psfrag{a}{{\small $A'$}} \psfrag{b}{{\small $B'$}}
    \psfrag{c}{{\small $C'$}} \psfrag{d}{{\small $D'$}}
    \centerline{\includegraphics[scale=0.7]{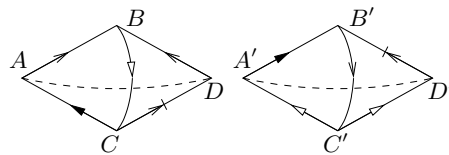}}
    \caption{An unallowable way of joining two tetrahedra along two faces}
    \label{fig-doubleedgebad}
    \end{figure}

    In this remaining possibility, face {\em ABC} is identified with
    $B'C'A'$ and face {\em BCD} is identified with $C'D'B'$, i.e., both
    identifications involve a rotation but in opposite directions.  In
    this case, faces {\em ABD} and $A'C'D'$ each have two edges
    identified to form a cone in contradiction to Lemma~\ref{l-conefaces}.
\end{proof}

\subsection{Properties of Face Pairing Graphs} \label{s-mainbadgraphs}

In contrast to Section~\ref{s-mainstruct} which analyses the structures
of triangulations based upon their face pairing graphs, in this
section we derive properties of the face pairing graphs themselves.
More specifically, we identify various types of subgraph that can never
occur within a face pairing graph of a closed minimal {\ppirr}
triangulation.

\begin{theorem} \label{t-censustripleedge}
    Let $G$ be a face pairing graph on $\geq 3$ vertices.
    If $G$ contains a {\em triple edge} (two vertices joined
    by three distinct edges as illustrated in
    Figure~\ref{fig-censustripleedge}), then $G$ cannot be the
    face pairing graph of a closed minimal {\ppirr} triangulation.

    \begin{figure}[htb]
    \centerline{\includegraphics[scale=0.7]{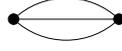}}
    \caption{A triple edge in a face pairing graph}
    \label{fig-censustripleedge}
    \end{figure}
\end{theorem}

\begin{proof}
    Suppose that $T$ is a closed minimal {\ppirr} triangulation with
    $\geq 3$ tetrahedra whose face pairing graph contains a triple edge.
    Observe that this triple edge corresponds to two distinct tetrahedra
    of $T$ that are joined along three different faces.
    We enumerate all possible ways in which this can be
    done and derive a contradiction in each case.

    Let these two tetrahedra be {\em ABCD} and $A'B'C'D'$ as illustrated
    in Figure~\ref{fig-censustripleedgetets}.  Faces {\em ABD} and
    $A'B'D'$ are identified, faces {\em BCD} and $B'C'D'$ are identified
    and faces {\em CAD} and $C'A'D'$ are identified, though not
    necessarily with these specific vertex identifications; faces may be
    rotated or reflected before they are identified.

    \begin{figure}[htb]
    \psfrag{A}{{\small $A$}} \psfrag{B}{{\small $B$}}
    \psfrag{C}{{\small $C$}} \psfrag{D}{{\small $D$}}
    \psfrag{a}{{\small $A'$}} \psfrag{b}{{\small $B'$}}
    \psfrag{c}{{\small $C'$}} \psfrag{d}{{\small $D'$}}
    \centerline{\includegraphics[scale=0.7]{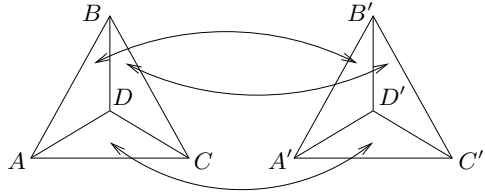}}
    \caption{Two tetrahedra to be joined along three faces}
    \label{fig-censustripleedgetets}
    \end{figure}

    Faces {\em ABC} and $A'B'C'$ remain unaccounted for, though they cannot
    be identified with each other since this would produce a
    2-tetrahedron triangulation.  Thus faces {\em ABC} and $A'B'C'$ are
    distinct faces of $T$; we refer to these as the {\em boundary
    faces} since they bound the subcomplex formed by the two tetrahedra
    under investigation.

    Each specific method of joining our two tetrahedra along the three
    pairs of faces is denoted by a {\em matching string}.  A matching
    string is a sequence of three symbols representing the
    transformations that are applied to faces {\em ABD}, {\em BCD} and
    {\em CAD} respectively before they are identified with their
    counterparts from the second tetrahedron.  Each symbol is one of the
    following.

    \begin{itemize}
        \item $\iota$: No transformation is applied.
        \item $\kappa$: The face is rotated clockwise.
        \item $\alpha$: The face is rotated anticlockwise.
        \item $c$: The face is reflected so that the centre point of the
        diagram (i.e., point $D$) remains fixed.
        \item $l$: The face is reflected so that the point at the clockwise
        end of the face (e.g., point $B$ on face {\em ABD}) remains fixed.
        \item $r$: The face is reflected so that the point at the
        anticlockwise end of the face (e.g., point $A$ on face {\em ABD})
        remains fixed.
    \end{itemize}

    A full list of precise face identifications corresponding to each symbol
    is given in Table~\ref{tab-censustriplefaceids}.

    \begin{table}[htb]
    \begin{center} $\begin{array}{|c|cccccc|}
        \hline
        & \iota & \kappa & \alpha & c & l & r \\
        \hline
        ABD & A'B'D' & B'D'A' & D'A'B' & B'A'D' & D'B'A' & A'D'B' \\
        BCD & B'C'D' & C'D'B' & D'B'C' & C'B'D' & D'C'B' & B'D'C' \\
        CAD & C'A'D' & A'D'C' & D'C'A' & A'C'D' & D'A'C' & C'D'A' \\
        \hline
    \end{array}$ \end{center}
    \caption{Precise face identifications corresponding to each
        transformation symbol}
    \label{tab-censustriplefaceids}
    \end{table}

    We can then enumerate the matching strings for all possible ways of
    joining our two tetrahedra along these three pairs of faces.  Each
    matching string is considered only once up to equivalence, where
    equivalence includes rotating the two tetrahedra, reflecting the two
    tetrahedra and swapping the two tetrahedra.

    As an example, a list of all matching strings up to equivalence for
    which the three face identifications have consistent orientations
    (i.e., which do not contain both a symbol
    from $\{\iota,\kappa,\alpha\}$ and a symbol from $\{c,l,r\}$) is as
    follows.
    \[
        \iota\iota\iota,~ \iota\iota\kappa,~ \iota\kappa\kappa,~
        \iota\kappa\alpha,~ \kappa\kappa\kappa,~ \kappa\kappa\alpha,~
        ccc,~ ccl,~ cll,~ clr,~ crl,~ lll,~ llr.
    \]
    Note that matching strings $\iota\kappa\alpha$ and
    $\iota\alpha\kappa$ are equivalent by
    swapping the two tetrahedra, but that $clr$ and $crl$ are not.

    A full list of all matching strings up to equivalence, allowing for both
    orientation-preserving and orientation-reversing face identifications,
    is too large to analyse by hand.  We can however use
    Lemma~\ref{l-censusdoubleedge} to restrict this list to a more
    manageable size.

    Specifically, Lemma~\ref{l-censusdoubleedge} imposes conditions upon
    how two tetrahedra may be joined along two faces.  In a matching
    string, each pair of adjacent symbols represents a joining of the
    two tetrahedra along two faces.  We can thus use this lemma to
    determine which symbols may be followed by which other symbols in
    a matching string.  The results of this analysis are presented in
    Table~\ref{tab-adjmatchingsymbols}.

    \begin{table}[htb]
    \begin{center} \begin{tabular}{|c|c|c|}
        \hline
        \bf Symbol & \bf May be followed by & \bf Cannot be followed by \\
        \hline
        $\iota$ & $\kappa$, $\alpha$ & $\iota$, $c$, $l$, $r$ \\
        $\kappa$ & $\iota$, $\kappa$, $c$, $r$ & $\alpha$, $l$ \\
        $\alpha$ & $\iota$, $\alpha$, $c$, $r$ & $\kappa$, $l$ \\
        $c$ & $\kappa$, $\alpha$, $l$, $r$ & $\iota$, $c$ \\
        $l$ & $\kappa$, $\alpha$, $c$, $l$ & $\iota$, $r$ \\
        $r$ & $c$, $r$ & $\iota$, $\kappa$, $\alpha$, $l$ \\
        \hline
    \end{tabular} \end{center}
    \caption{Possible adjacent symbols within a matching string}
    \label{tab-adjmatchingsymbols}
    \end{table}

    Up to rotation of the two tetrahedra, this reduces our list of
    possible matching strings to the following.
    \[
        \iota\kappa\kappa,~\iota\alpha\alpha,~
        \kappa\kappa\kappa,~\alpha\alpha\alpha,~\kappa\kappa c,~
        \alpha\alpha c,~
        \kappa c l,~\kappa r c,~\alpha c l,~\alpha r c,~
        cll,~crr,~lll,~rrr.
    \]
    Allowing for the two tetrahedra to be reflected and/or swapped, this
    list can be further reduced to just the following six matching strings.
    \[
        \iota\kappa\kappa,~\kappa\kappa\kappa,~\kappa\kappa c,~
        \kappa c l,~cll,~lll.
    \]

    These six matching strings can be split into three
    categories, where the matchings in each category give rise to a
    similar contradiction using almost identical arguments.  We
    examine each category in turn.

    \begin{itemize}
        \item {\bf Conical faces ($\iota\kappa\kappa$, $\kappa c l$):}
        Simply by following the edge identifications induced by our
        chosen face identifications, we can see that each of these
        matchings gives rise to a face containing two edges that are
        identified to form a cone.  Lemma~\ref{l-conefaces} then
        provides us with our contradiction.

        An example of this is illustrated in
        Figure~\ref{fig-censustripleedgekcl} which shows the induced edge
        identifications for the matching $\kappa cl$.  Conical
        faces in this example include {\em ABD} and the boundary face
        {\em ABC}.

        \begin{figure}[htb]
        \psfrag{A}{{\small $A$}} \psfrag{B}{{\small $B$}}
        \psfrag{C}{{\small $C$}} \psfrag{D}{{\small $D$}}
        \psfrag{a}{{\small $A'$}} \psfrag{b}{{\small $B'$}}
        \psfrag{c}{{\small $C'$}} \psfrag{d}{{\small $D'$}}
        \centerline{\includegraphics[scale=0.7]{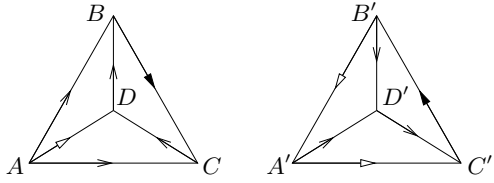}}
        \caption{Edge identifications for the matching $\kappa cl$}
        \label{fig-censustripleedgekcl}
        \end{figure}

        \item {\bf Spherical subcomplexes ($\kappa\kappa\kappa$,
        $\kappa\kappa c$):}
        Again we follow the induced edge identifications, but in this
        case the consequence is that the two boundary faces {\em ABC}
        and $A'B'C'$ are joined at their edges to form a two-triangle sphere
        in contradiction to Lemma~\ref{l-splittripillow}.
        This behaviour is illustrated in Figure~\ref{fig-censustripleedgekkc}
        for the matching $\kappa\kappa c$.

        \begin{figure}[htb]
        \psfrag{A}{{\small $A$}} \psfrag{B}{{\small $B$}}
        \psfrag{C}{{\small $C$}} \psfrag{D}{{\small $D$}}
        \psfrag{a}{{\small $A'$}} \psfrag{b}{{\small $B'$}}
        \psfrag{c}{{\small $C'$}} \psfrag{d}{{\small $D'$}}
        \centerline{\includegraphics[scale=0.7]{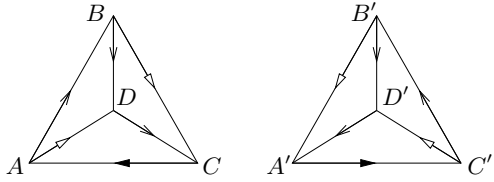}}
        \caption{Edge identifications for the matching $\kappa\kappa c$}
        \label{fig-censustripleedgekkc}
        \end{figure}

        \item {\bf Bad vertex links ($cll$, $lll$):}
        For these matching strings, all eight vertices $A$,
        $B$, $C$, $D$, $A'$, $B'$, $C'$ and $D'$ are identified as a
        single vertex in the triangulation.  This vertex
        lies on the boundary faces {\em ABC} and $A'B'C'$, and so the vertex
        link as restricted to our two tetrahedra will be incomplete
        (i.e., will have boundary components).  However, for
        each of these matchings, this
        partial vertex link is observed to be a once-punctured
        torus.  Therefore, however the entire triangulation $T$ is formed,
        the complete vertex link in $T$ cannot be a sphere (since
        there is no way to fill in the boundary of a
        punctured torus to form a sphere).  Thus again we have a
        contradiction since $T$ cannot be a triangulation of a
        closed 3-manifold.

        As an example of this behaviour,
        Figure~\ref{fig-censustripleedgelll} illustrates the induced
        edge and vertex identifications for matching $lll$ as well as
        the corresponding vertex link.  The link is shown as eight
        individual triangles followed by a combined figure, which we see
        is indeed a punctured torus.

        \begin{figure}[htb]
        \psfrag{A}{{\small $A$}} \psfrag{B}{{\small $B$}}
        \psfrag{C}{{\small $C$}} \psfrag{D}{{\small $D$}}
        \psfrag{a}{{\small $A'$}} \psfrag{b}{{\small $B'$}}
        \psfrag{c}{{\small $C'$}} \psfrag{d}{{\small $D'$}}
        \centerline{\includegraphics[scale=0.7]{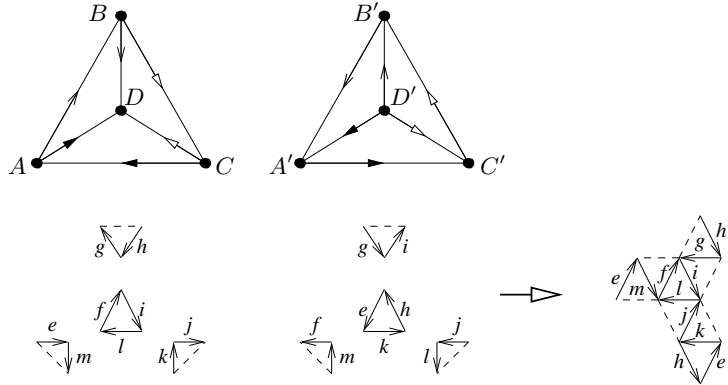}}
        \caption{The partial vertex link for the matching $lll$}
        \label{fig-censustripleedgelll}
        \end{figure}
    \end{itemize}

    Thus we see that every method of identifying three faces of the
    first tetrahedron with three faces of the second gives rise to a
    contradiction, and so our result is established.
\end{proof}

\begin{theorem} \label{t-censusbrokenchain}
    Let $G$ be a face pairing graph on $\geq 3$ vertices.
    If $G$ contains as a subgraph a {\em broken double-ended chain}
    (a double-ended chain missing one interior edge as illustrated in
    Figure~\ref{fig-censusbrokenchain}) and if $G$ is
    not simply a double-ended chain itself, then $G$ cannot be the
    face pairing graph of a closed minimal {\ppirr} triangulation.

    \begin{figure}[htb]
    \centerline{\includegraphics[scale=0.7]{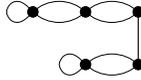}}
    \caption{A broken double-ended chain in a face pairing graph}
    \label{fig-censusbrokenchain}
    \end{figure}
\end{theorem}

\begin{proof}
    Observe that a broken double-ended chain is merely a pair of
    one-ended chains
    joined by an edge.  Let $T$ be a closed minimal {\ppirr}
    triangulation whose face pairing graph contains a broken
    double-ended chain.  Then
    Theorem~\ref{t-censuschainoneend} implies that $T$ contains a pair of
    layered solid tori whose boundaries are joined along one face.

    This situation is depicted in the left hand diagram of
    Figure~\ref{fig-censuslstlstjoin}, where the two torus boundaries
    of the layered solid tori are shown and where face {\em ABC} is
    identified with face {\em XYZ}.  The resulting edge identifications of
    the remaining two boundary faces are illustrated in the right hand
    diagram of this figure; in particular, it can be seen that these
    remaining boundary faces form a two-triangle sphere (though with all
    three vertices pinched together, since each layered solid torus has
    only one vertex).

    \begin{figure}[htb]
    \psfrag{A}{{\small $A$}} \psfrag{B}{{\small $B$}}
    \psfrag{C}{{\small $C$}}
    \psfrag{X}{{\small $X$}} \psfrag{Y}{{\small $Y$}}
    \psfrag{Z}{{\small $Z$}}
    \psfrag{F1}{{\small $F_1$}} \psfrag{F2}{{\small $F_2$}}
    \centerline{\includegraphics[scale=0.7]{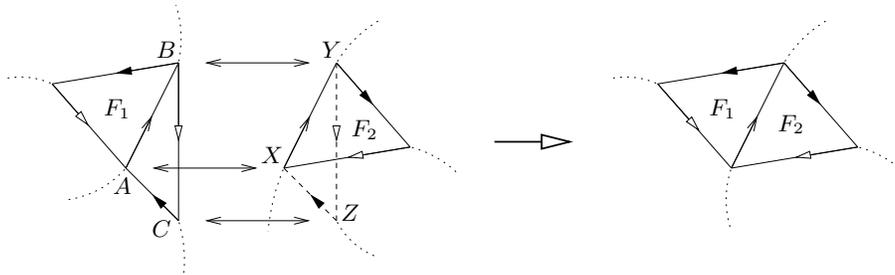}}
    \caption{Two layered solid tori joined along a face}
    \label{fig-censuslstlstjoin}
    \end{figure}

    Let these remaining boundary faces be $F_1$ and $F_2$.  If these
    faces are not identified then $F_1$ and $F_2$ satisfy the conditions
    of Lemma~\ref{l-splittripillow} and so $T$ cannot be both minimal
    and \ppirr.  Therefore faces $F_1$ and $F_2$ are identified.
    Returning to the face pairing graph for $T$, this implies that the
    single edge between the two one-ended chains is in fact a double
    edge and so the graph contains an entire double-ended chain.

    However, since every vertex in a double-ended chain has degree 4,
    Lemma~\ref{l-graphprops} shows this double-ended chain to be
    the entire face pairing graph, contradicting the initial conditions of
    this theorem.
\end{proof}

\begin{theorem} \label{t-censuschaindoublehandle}
    Let $G$ be a face pairing graph on $\geq 3$ vertices.
    If $G$ contains as a subgraph a {\em one-ended chain with a double
    handle} (a double-ended chain with one end edge replaced by a triangle
    containing one double edge as illustrated in
    Figure~\ref{fig-censuschaindoublehandle}), then $G$ cannot be the
    face pairing graph of a closed minimal {\ppirr} triangulation.

    \begin{figure}[htb]
    \centerline{\includegraphics[scale=0.7]{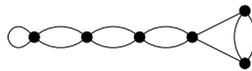}}
    \caption{A one-ended chain with a double handle in a face pairing graph}
    \label{fig-censuschaindoublehandle}
    \end{figure}
\end{theorem}

\begin{proof}
    Let $T$ be a closed minimal {\ppirr} triangulation with
    $\geq 3$ tetrahedra whose face pairing graph contains a one-ended
    chain with a double handle.  From Theorem~\ref{t-censuschainoneend} we
    see that the one-ended chain must correspond to a layered solid torus in
    $T$.  The double handle in turn must correspond to two additional
    tetrahedra each of which is joined to the other along two faces and
    each of which is joined to one of the boundary faces of the layered
    solid torus.

    \begin{figure}[htb]
    \centerline{\includegraphics[scale=0.7]{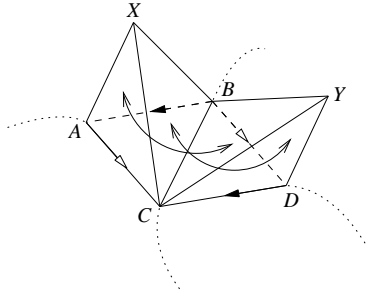}}
    \caption{Tetrahedra corresponding to a one-ended chain with a
        double handle}
    \label{fig-censusdoublehandletri}
    \end{figure}

    This construction is illustrated in
    Figure~\ref{fig-censusdoublehandletri}.
    The layered solid torus lies beneath faces {\em ABC} and {\em BCD} which
    form its torus boundary.  The two additional
    tetrahedra are {\em XABC} and {\em YBCD}; observe that each of these
    tetrahedra is joined to one of the boundary faces of the layered
    solid torus.  The two new tetrahedra are then joined to each other
    along two faces; in this example faces {\em XAC} and {\em YCB} are
    identified and faces {\em XCB} and {\em DCY} are identified, though
    different pairs of faces may be used.

    As in the proof of Theorem~\ref{t-censustripleedge}, we
    enumerate all possible ways in which this construction can
    be carried out and in each case derive a contradiction.  To
    assist in our task we describe a simple way of representing
    each possible variant of this construction.

    The scenario presented in Figure~\ref{fig-censusdoublehandletri}
    can be distilled into a simplified diagram as illustrated in
    Figure~\ref{fig-censusdoublehandleeg}.  We begin with the
    two-triangle torus that forms the boundary of the layered solid
    torus as shown in the left hand diagram of
    Figure~\ref{fig-censusdoublehandleeg}.  We then add our two new
    tetrahedra, converting the two-triangle torus into a six-triangle
    torus as illustrated in the central diagram of
    Figure~\ref{fig-censusdoublehandleeg}.  Finally we add markings to
    the diagram to illustrate how the two new tetrahedra are to be
    joined along two faces.
    This is done by marking vertices $a$, $b$ and $c$ of the
    first face and vertices $x$, $y$ and $z$ of the second face.
    It can be seen from the right
    hand diagram of Figure~\ref{fig-censusdoublehandleeg} that
    faces {\em XAC} and {\em YCB} are identified and faces
    {\em XCB} and {\em DCY} are identified as described earlier.
    For clarity the two faces that have not yet been identified with any
    others are shaded.

    \begin{figure}[htb]
    \centerline{\includegraphics[scale=0.5]{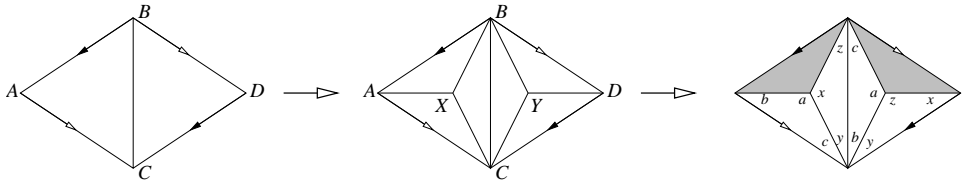}}
    \caption{A simplified diagram for a one-ended chain with a double handle}
    \label{fig-censusdoublehandleeg}
    \end{figure}

    In a similar fashion we can represent any set of tetrahedra
    corresponding to a one-ended chain with a double handle
    using a diagram similar to the right hand diagram of
    Figure~\ref{fig-censusdoublehandleeg}.  In this way we can
    enumerate all possible diagrams and in each case prove that $T$
    cannot be a closed minimal {\ppirr} triangulation.

    Recall from Lemma~\ref{l-censusdoubleedge} that if $T$ is a closed
    minimal {\ppirr} triangulation then there are restrictions upon the
    possible ways in which our two new tetrahedra can be joined along
    two faces.  Furthermore, Lemma~\ref{l-tetringor} requires that if
    two adjacent faces are to be identified then this must be done in an
    orientation-preserving manner.  By
    ignoring all diagrams that do not conform to
    Lemmas~\ref{l-tetringor} and~\ref{l-censusdoubleedge} and by
    exploiting the symmetries of the six-triangle torus and the layered
    solid torus that lies beneath it, we can reduce the set of all possible
    diagrams to those depicted in
    Figures~\ref{fig-censusdoublehandlecasesor}
    and~\ref{fig-censusdoublehandlecasesnor}.
    Figure~\ref{fig-censusdoublehandlecasesor} contains the diagrams
    for constructions that are orientation-preserving and
    Figure~\ref{fig-censusdoublehandlecasesnor} contains the diagrams
    corresponding to non-orientable structures.

    \begin{figure}[htb]
    \psfrag{A1}{$\alpha_1$} \psfrag{A2}{$\alpha_2$} \psfrag{A3}{$\alpha_3$}
    \psfrag{B1}{$\beta_1$} \psfrag{B2}{$\beta_2$} \psfrag{B3}{$\beta_3$}
    \psfrag{C1}{$\gamma_1$} \psfrag{C2}{$\gamma_2$} \psfrag{C3}{$\gamma_3$}
    \psfrag{D1}{$\delta_1$} \psfrag{D2}{$\delta_2$}
    \centerline{\includegraphics[scale=0.5]{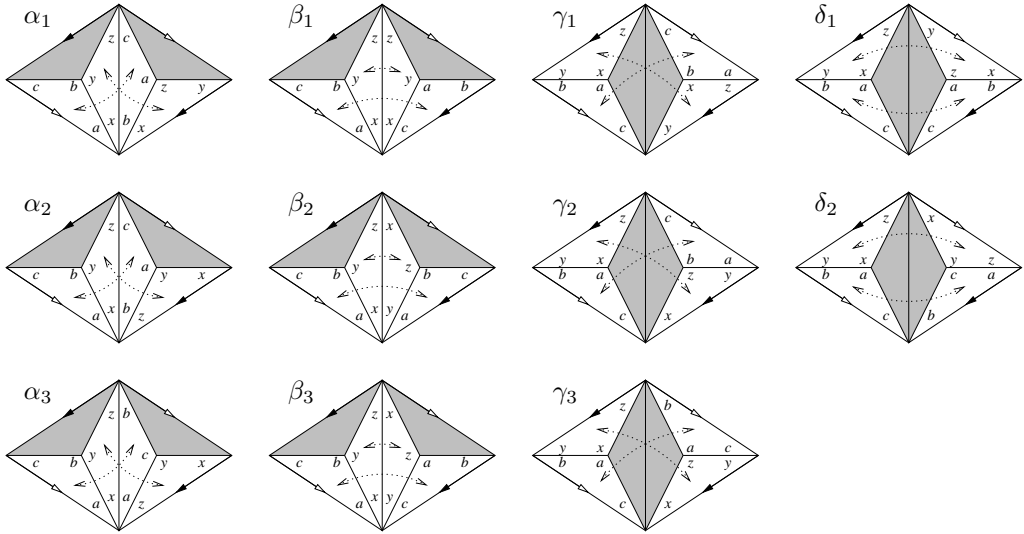}}
    \caption{All orientable diagrams for a one-ended chain with a double handle}
    \label{fig-censusdoublehandlecasesor}
    \end{figure}

    \begin{figure}[htb]
    \psfrag{A4}{$\alpha_4$} \psfrag{A5}{$\alpha_5$}
    \psfrag{A6}{$\alpha_6$} \psfrag{A7}{$\alpha_7$}
    \psfrag{B4}{$\beta_4$} \psfrag{B5}{$\beta_5$}
    \psfrag{D3}{$\delta_3$} \psfrag{D4}{$\delta_4$}
    \centerline{\includegraphics[scale=0.5]{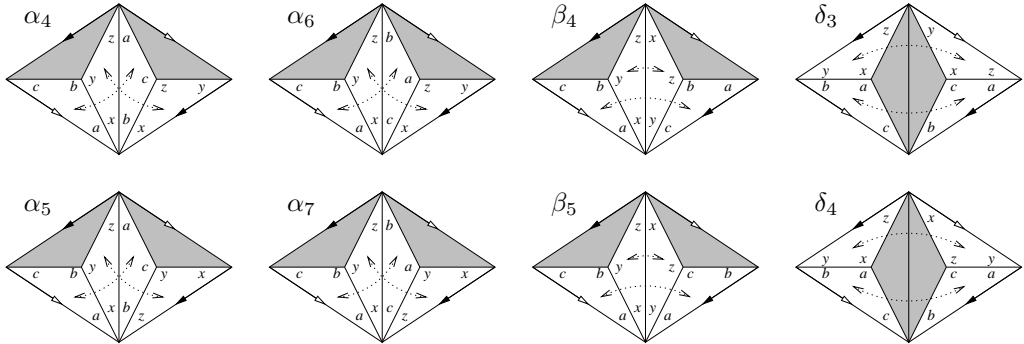}}
    \caption{All non-orientable diagrams for a one-ended chain with a
        double handle}
    \label{fig-censusdoublehandlecasesnor}
    \end{figure}

    As in the proof of Theorem~\ref{t-censustripleedge} we can divide
    our 19 different diagrams $\alpha_1,\ldots,\delta_4$ into a
    small number of categories, where the diagrams in each category give
    rise to similar contradictions using almost identical arguments.
    To assist with this process the edge identifications induced in each
    diagram by the corresponding face identifications are shown in
    Figures~\ref{fig-censusdoublehandlearrowsor}
    and~\ref{fig-censusdoublehandlearrowsnor}.
    The diagrams can then be split into categories as follows.

    \begin{figure}[htb]
    \psfrag{A1}{$\alpha_1$} \psfrag{A2}{$\alpha_2$} \psfrag{A3}{$\alpha_3$}
    \psfrag{B1}{$\beta_1$} \psfrag{B2}{$\beta_2$} \psfrag{B3}{$\beta_3$}
    \psfrag{C1}{$\gamma_1$} \psfrag{C2}{$\gamma_2$} \psfrag{C3}{$\gamma_3$}
    \psfrag{D1}{$\delta_1$} \psfrag{D2}{$\delta_2$}
    \centerline{\includegraphics[scale=0.5]{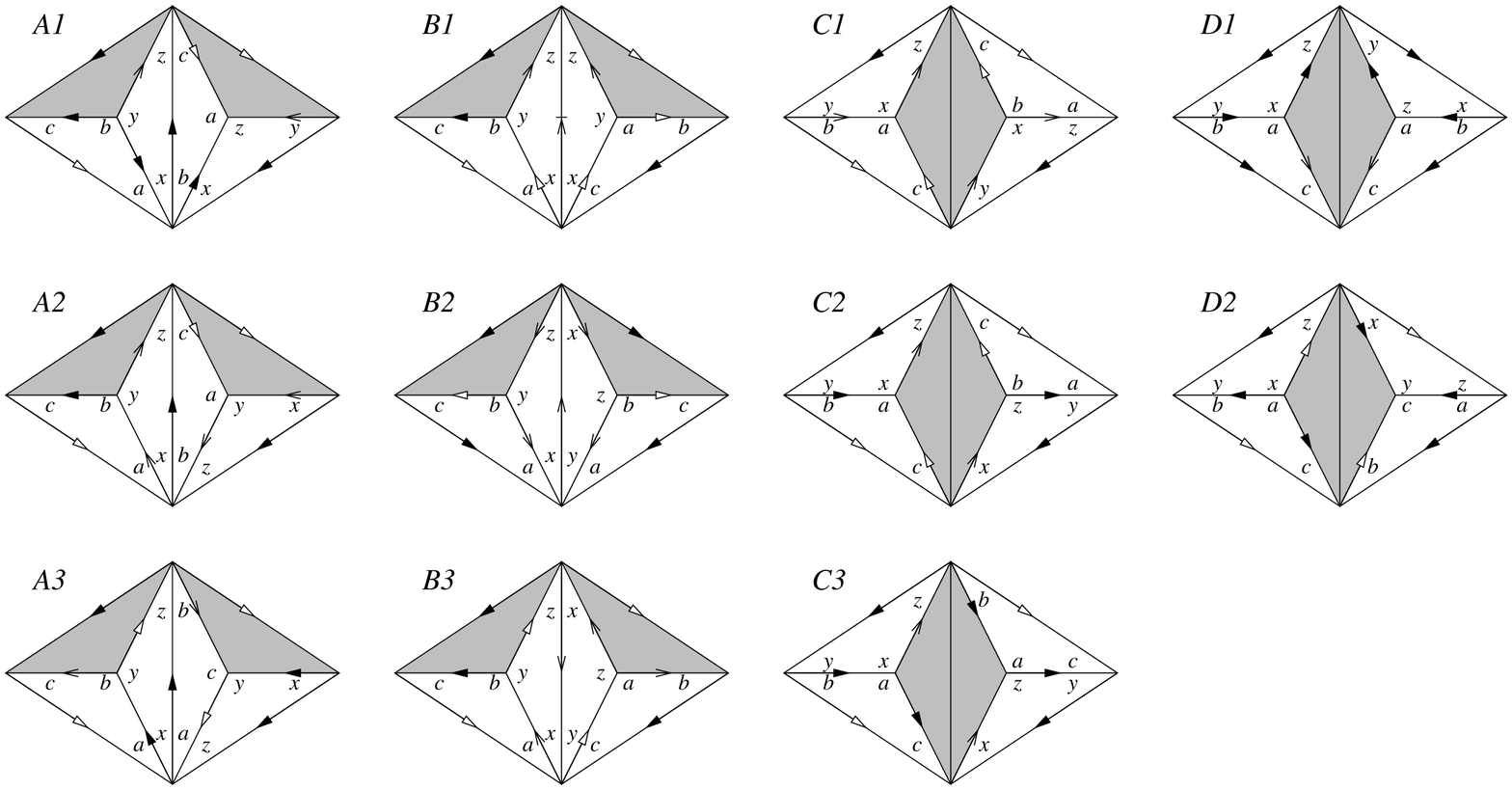}}
    \caption{Edge identifications for the orientable face
        identification diagrams}
    \label{fig-censusdoublehandlearrowsor}
    \end{figure}

    \begin{figure}[htb]
    \psfrag{A4}{$\alpha_4$} \psfrag{A5}{$\alpha_5$}
    \psfrag{A6}{$\alpha_6$} \psfrag{A7}{$\alpha_7$}
    \psfrag{B4}{$\beta_4$} \psfrag{B5}{$\beta_5$}
    \psfrag{D3}{$\delta_3$} \psfrag{D4}{$\delta_4$}
    \centerline{\includegraphics[scale=0.5]{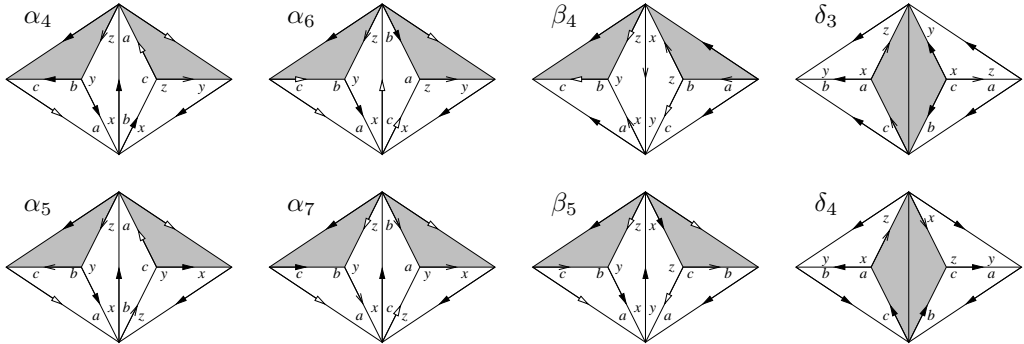}}
    \caption{Edge identifications for the non-orientable face
        identification diagrams}
    \label{fig-censusdoublehandlearrowsnor}
    \end{figure}

    \begin{itemize}
        \item {\bf Conical faces ($\alpha_1$, $\alpha_2$,
            $\alpha_4$, $\alpha_6$, $\beta_1$, $\beta_3$, $\beta_4$,
            $\delta_3$, $\delta_4$):}
        In each of these diagrams we find at least one face with two
        edges identified to form a cone.  From Lemma~\ref{l-conefaces}
        it follows that $T$ cannot be a closed minimal
        {\ppirr} triangulation.

        \item {\bf Spherical subcomplexes ($\alpha_3$, $\beta_2$,
            $\beta_5$, $\delta_1$):}
        In each of these diagrams we can observe that the two shaded
        faces are joined at their edges to form a two-triangle sphere as
        described in Lemma~\ref{l-splittripillow}, again contradicting
        the claim that $T$ is a closed
        minimal {\ppirr} triangulation.

        \item {\bf Bad vertex links ($\alpha_5$, $\alpha_7$,
            $\gamma_1$, $\gamma_2$, $\gamma_3$, $\delta_2$):}
        In each of these diagrams it can be observed that all six
        vertices illustrated are in fact identified
        as a single vertex in the triangulation.  We can calculate the
        link of this vertex as restricted to the portion of the
        triangulation that we are examining.  For diagrams
        $\gamma_1$, $\gamma_2$ and $\gamma_3$ this link is
        observed to be a once-punctured torus, for diagram
        $\delta_2$ it is observed to be a once-punctured genus two torus
        and for diagrams $\alpha_5$ and $\alpha_7$ it is observed to be
        non-orientable.
        In each of these cases, however the
        remainder of the triangulation is constructed it is impossible
        for this vertex link to be extended to become a sphere.
        Thus $T$ cannot be a triangulation of a closed 3-manifold.

        The vertex link calculation is illustrated in
        Figure~\ref{fig-censusdoublehandlelink} for diagram $\gamma_1$.
        The disc on the left with edges $p$, $q$, $r$ and $s$
        represents the vertex link of the layered
        solid torus and the triangles beside it represent the pieces of
        vertex link taken from the two new tetrahedra.  These pieces
        are combined into a single surface on the right hand side of the
        diagram which we see is indeed a once-punctured torus.

        \begin{figure}[htb]
        \centerline{\includegraphics[scale=0.7]{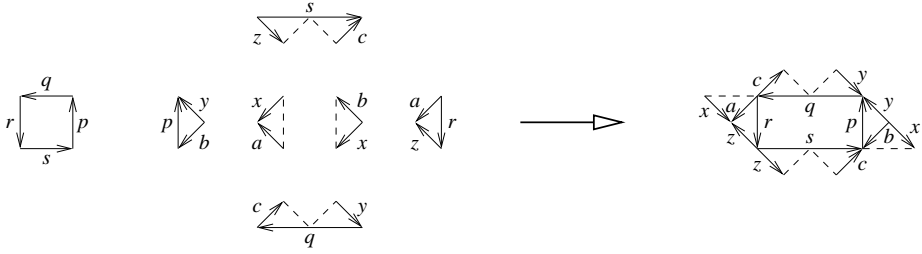}}
        \caption{Calculating the vertex link for diagram $\gamma_1$}
        \label{fig-censusdoublehandlelink}
        \end{figure}
    \end{itemize}

    Thus for each diagram of Figure~\ref{fig-censusdoublehandlecasesor}
    and Figure~\ref{fig-censusdoublehandlecasesnor} we
    observe that $T$ cannot be a closed minimal {\ppirr} triangulation.
\end{proof}

\section{Enumeration of 3-Manifold Triangulations} \label{s-enumeration}

We move now to the task of applying the results of Section~\ref{s-main} to
the enumeration of closed minimal {\ppirr} triangulations.  As described
in the introduction, the enumeration of 3-manifold triangulations is a
slow procedure typically based upon a brute force search through
possible identifications of tetrahedron faces.  Our aim then is to use
the face pairing graph results of Section~\ref{s-main} to impose
restrictions upon this brute force search and thus improve the
performance of the algorithm.

Table~\ref{tab-oldtimes} illustrates the performance of the enumeration
algorithm before these new improvements are introduced.
Included in the table are the counts up to isomorphism of
closed minimal {\ppirr} triangulations formed
from various numbers of tetrahedra, as well as the
time taken by the computer to enumerate these triangulations.  All times
are displayed as {\em h}:{\em mm}:{\em ss} and are measured on a single
1.2GHz Pentium~III processor.  A
time of 0:00 simply indicates a running time of less than half a second.

\begin{table}[htb]
\begin{center} \begin{tabular}{|c|r|r|r|r|}
    \hline
    & \multicolumn{2}{|c|}{\bf Orientable} &
        \multicolumn{2}{|c|}{\bf Non-Orientable} \\
    \hline
    \bf Tetrahedra &
        \bf Triang. & \bf Running Time &
        \bf Triang. & \bf Running Time \\
    \hline
    3 & 7 & 0:00 & 0 & 0:02 \\
    4 & 15 & 0:03 & 0 & 3:03 \\
    5 & 40 & 2:28 & 0 & 6:01:53 \\
    6 & 115 & 2:49:29 & 24 & 5 weeks \\
    \hline
\end{tabular} \end{center}
\caption{Enumeration statistics without face pairing graph improvements}
\label{tab-oldtimes}
\end{table}

The enumerations described in Table~\ref{tab-oldtimes} were performed by
{\regina} \cite{regina}, a freely available computer program that
can perform a variety of calculations and procedures in 3-manifold
topology.  Statistics are not offered for one or two tetrahedra since
these running times all round to 0:00.

In Section~\ref{s-algmeliminatepairings} we outline the structure of a
typical enumeration algorithm and use the graph constraints derived in
Section~\ref{s-mainbadgraphs} to make some immediate
efficiency improvements.  Section~\ref{s-algmredesign} then uses the
triangulation constraints derived in
Section~\ref{s-mainstruct} to redesign the enumeration algorithm for
a more significant increase in performance.
Finally in Section~\ref{s-statistics} we
present additional statistics to illustrate how well the newly designed
algorithm performs in practice.

\subsection{Eliminating Face Pairings} \label{s-algmeliminatepairings}

The algorithm for enumerating 3-manifold triangulations can be split
into two largely independent tasks, these being the generation of face
pairing graphs and the selection of rotations and reflections for the
corresponding face identifications.  This is a fairly
natural way of approaching an enumeration of triangulations and is seen
back in the earliest hyperbolic census of Hildebrand and Weeks
\cite{cuspedcensusold}.  Indeed this splitting of tasks appears
to be used by all enumeration algorithms described in the literature,
although different authors employ different variants of the general
technique.

More specifically, a blueprint for an enumeration algorithm can
be described as follows.

\begin{algorithm}[Enumeration of Triangulations] \label{a-census}
    To enumerate up to isomorphism all 3-manifold triangulations formed
    from $n$ tetrahedra that satisfy some particular census
    constraints (such as minimality and {\ppirrty}), we perform the
    following steps.

    \begin{enumerate}
        \item \label{en-census-pairings}
        Generate up to isomorphism all face pairings graphs for $n$
        tetrahedra, i.e., all connected 4-valent multigraphs on $n$
        vertices as described by Lemma~\ref{l-graphprops}.

        \item \label{en-census-gluings}
        For each generated face pairing graph, recursively try all
        possible rotations and reflections
        for each pair of faces to be identified.

        In general each pair of faces can be identified according to
        one of six possible rotations or reflections for non-orientable
        triangulations, or one of three possible rotations or
        reflections for orientable triangulations (allowing only the
        orientation-preserving identifications).

        Thus, with $2n$ pairs of identified faces, this requires a
        search through approximately $6^{2n}$ or $3^{2n}$ possible
        combinations of rotations and reflections for non-orientable
        or orientable triangulations respectively.  In practice this
        search can be pruned somewhat.  For instance, the results of
        Section~\ref{s-lowdegedges} show that the search tree can be
        pruned where it becomes apparent that an edge of low degree
        will be present in the final triangulation.  Other pruning
        techniques are described in the literature and frequently depend
        upon the specific census constraints under consideration.

        \item \label{en-census-verify}
        For each triangulation thus constructed, test whether it
        satisfies the full set of census constraints and if so then
        include it in the final list of results.
    \end{enumerate}
\end{algorithm}

In practise, the time spent generating face pairing graphs
(step~\ref{en-census-pairings}) is utterly negligible.  For the six
tetrahedron closed non-orientable census that consumed five weeks of
processor time as seen in Table~\ref{tab-oldtimes}, the generation of
face pairing graphs took under a hundredth of a second.

Our first improvement to Algorithm~\ref{a-census} is then as follows.
When generating face pairing graphs in step~\ref{en-census-pairings} of
Algorithm~\ref{a-census}, we throw away any face pairing graphs that do
not satisfy the properties proven in Section~\ref{s-mainbadgraphs}.
Specifically we throw away face pairing graphs containing triple edges
(Theorem~\ref{t-censustripleedge}), broken double-ended chains
(Theorem~\ref{t-censusbrokenchain}) or one-ended
chains with double handles (Theorem~\ref{t-censuschaindoublehandle}).

It is worth noting that Martelli and Petronio \cite{italian9} employ a
technique that likewise eliminates face pairing graphs for the purpose
of enumerating bricks (a particular class of subcomplexes within
3-manifold triangulations).  Their results are phrased in terms of
special spines of 3-manifolds and relate to orientable triangulations
satisfying specific structural constraints.

Since the enumeration of face pairing graphs is so fast, it seems plausible
that if we can eliminate $p\%$ of all face pairing graphs in this way
then we can eliminate approximately $p\%$ of the total running time
of the algorithm.  Table~\ref{tab-badfacepairings} illustrates how
effective this technique is in practice.  Specifically, it lists how
many face pairing graphs on $n$ vertices for $3 \leq n \leq 11$ contain
each of the unallowable structures listed above.  All counts are given
up to isomorphism and were obtained using the program {\regina}
\cite{regina}.

The individual columns of Table~\ref{tab-badfacepairings} have the
following meanings.
\begin{itemize}
    \item {\em Vertices:} The number of vertices $n$.
    \item {\em Total:} The total number $t$ of connected 4-valent
    multigraphs on $n$ vertices.
    \item {\em None:} The number and percentage of graphs from this
    total $t$ in which none of the undesirable structures listed above were
    found.
    \item {\em Some:} The number of graphs from this total $t$ in which at
    least one of these undesirable structures was found.
    \item {\em Triple:} The number of graphs from the total $t$ that
    contain a triple edge.
    \item {\em Broken:} The number of graphs from the total $t$ that
    contain a broken double-ended chain.
    \item {\em Handle:} The number of graphs from the total $t$ that
    contain a one-ended chain with a double handle.
    \item {\em Time:} The running time taken to calculate the values in
    this row of the table.  Running times are measured on a single
    1.2GHz Pentium~III processor and are displayed as
    {\em h}:{\em mm}:{\em ss}.
\end{itemize}

Note that the {\em Total} column should equal the {\em None} column plus
the {\em Some} column in each row.  Note also that the sum
of the {\em Triple}, {\em Broken} and {\em Handle} columns might exceed the
{\em Some} column since some graphs may contain more than one type of
undesirable structure.

\begin{table}[htb]
\begin{center} \begin{tabular}{|r|r|rr|r|r|r|r|r|}
    \hline
    \bf Vertices & \bf Total & \multicolumn{2}{|c|}{\bf None} & \bf Some &
        \bf Triple & \bf Broken & \bf Handle &
        \bf Time \\
    \hline
    3 & 4 & 2 & (50\%) & 2 & 1 & 1 & 1 & 0:00 \\
    4 & 10 & 4 & (40\%) & 6 & 3 & 3 & 2 & 0:00 \\
    5 & 28 & 12 & (43\%) & 16 & 8 & 10 & 4 & 0:00 \\
    6 & 97 & 39 & (40\%) & 58 & 29 & 36 & 12 & 0:00 \\
    7 & 359 & 138 & (38\%) & 221 & 109 & 137 & 40 & 0:01 \\
    8 & 1\,635 & 638 & (39\%) & 997 & 497 & 608 & 155 & 0:05 \\
    9 & 8\,296 & 3\,366 & (41\%) & 4\,930 & 2\,479 & 2\,976 & 685 & 0:44 \\
    10 & 48\,432 & 20\,751 & (43\%) & 27\,681 & 14\,101 & 16\,568 & 3\,396 &
        7:21 \\
    11 & 316\,520 & 143\,829 & (45\%) & 172\,691 & 88\,662  & 102\,498 &
        18\,974 & 1:20:48 \\
    \hline
\end{tabular} \end{center}
\caption{Frequency of undesirable structures within face pairing graphs}
\label{tab-badfacepairings}
\end{table}

Examining Table~\ref{tab-badfacepairings} we see then that for each
number of vertices a little over half of the possible face pairing
graphs can be eliminated using the results of
Section~\ref{s-mainbadgraphs}.  Whilst this does not
reduce the time complexity of the algorithm, it does allow us
to reduce the running time by more than 50\%, a significant improvement
for a census that may take months or years to complete.

\begin{figure}[htb]
\centerline{\includegraphics[scale=0.7]{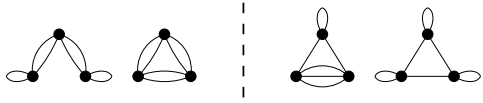}}
\caption{All possible face pairing graphs on three vertices}
\label{fig-censusfacegraphs3}
\end{figure}

\begin{figure}[htb]
\centerline{\includegraphics[scale=0.7]{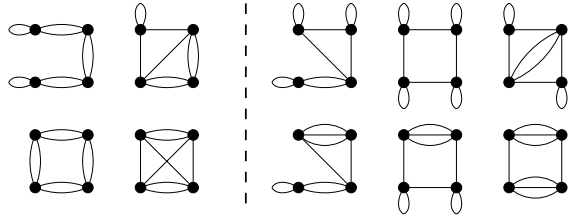}}
\caption{All possible face pairing graphs on four vertices}
\label{fig-censusfacegraphs4}
\end{figure}

It is worth examining how comprehensive the theorems of
Section~\ref{s-mainbadgraphs} are.  We turn our attention here
specifically to the orientable case.  Figure~\ref{fig-censusfacegraphs3}
lists all possible 3-vertex face pairing graphs, i.e., all connected
4-valent graphs on three vertices.  A list of all possible 4-vertex face
pairing graphs is provided in Figure~\ref{fig-censusfacegraphs4}.  In
each figure, the graphs to the left of the dotted line all lead to
closed orientable minimal {\ppirr} triangulations whereas the graphs to
the right of the dotted line do not.  It can be seen that in both figures,
every right hand graph contains either a triple edge, a broken
double-ended chain or a one-ended chain with a double handle.

Thus for three and four tetrahedra, the results of
Section~\ref{s-mainbadgraphs} in fact perfectly divide the face
pairings into those that lead to desirable triangulations and those that
do not.

\begin{figure}[htb]
\centerline{\includegraphics[scale=0.7]{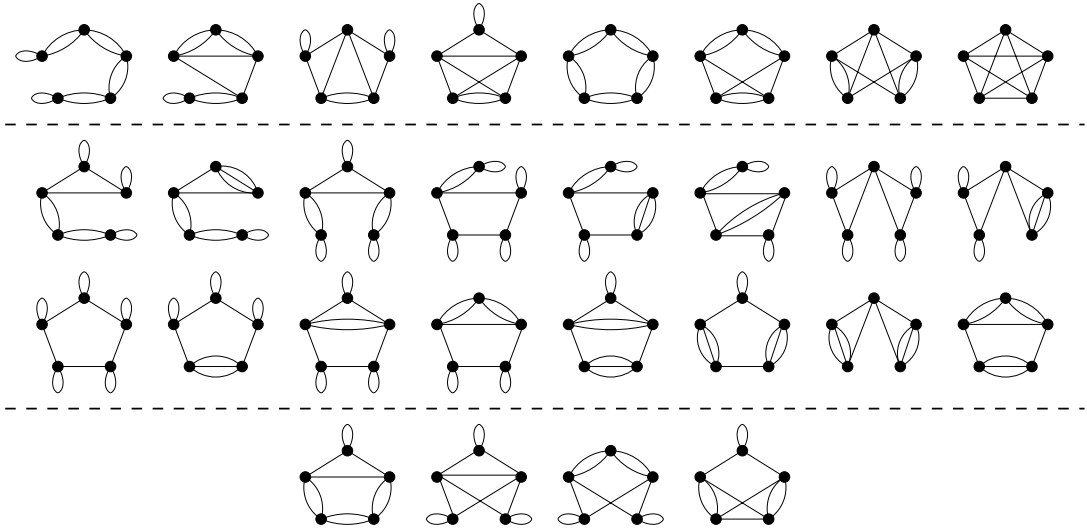}}
\caption{All possible face pairing graphs on five vertices}
\label{fig-censusfacegraphs5}
\end{figure}

For five tetrahedra these theorems no longer perfectly divide the face
pairings as we would like.  Figure~\ref{fig-censusfacegraphs5} lists all
5-vertex face pairing graphs; there are 28 in total.  The 8 graphs in
the top section all lead to closed orientable minimal {\ppirr}
triangulations and the remaining 20 do not.  Of these remaining 20, the
16 graphs in the middle section each contain a triple edge, a broken
double-ended chain or a one-ended chain with a double handle.  We see
then that the final 4 graphs in the bottom section can never lead to a
desirable triangulation but are not identified as such by the theorems
of Section~\ref{s-mainbadgraphs}.  Thus there remains
more work to be done in the analysis of face pairing graphs.

\subsection{Improved Generation of Rotations and Reflections}
    \label{s-algmredesign}

Consider once more Algorithm~\ref{a-census}.  This describes the enumeration
algorithm as a three-stage process in which we generate face pairing graphs,
generate rotations and reflections for each face pairing graph and then
analyse the resulting triangulations.
Recall from Section~\ref{s-algmeliminatepairings} that the running time spent
generating face pairing graphs is inconsequential in the
context of the overall enumeration algorithm.  In fact almost the entire
running time of the algorithm is spent in the second stage generating
rotations and reflections.

As seen in step~\ref{en-census-gluings} of Algorithm~\ref{a-census} this
generation of rotations and reflections is a recursive process in which
we choose
from six possibilities (or three for an orientable triangulation) for
each pair of faces that are to be identified.  This allows for up to
$6^{2n}$ (or $3^{2n}$ for an orientable triangulation) different
combinations of rotations and reflections
for any given face pairing graph.

It is thus in the generation of rotations and reflections that we
should seek to make the strongest improvements to the enumeration algorithm.
Algorithm~\ref{a-census} already suggests techniques
through which this recursion can be pruned and thus made faster.
It is possible however that instead of chipping away at the recursion
with more and more intricate pruning techniques, we could perhaps
make more significant improvements by substantially redesigning
the recursion using the face pairing graph results of
Section~\ref{s-mainstruct}.

Recall from Theorem~\ref{t-censuschainoneend} that every one-ended chain
in a face pairing graph must correspond to a layered solid torus in the
resulting triangulation.  Recall also from
Lemma~\ref{l-censusdoubleedge} that every double edge in a face pairing
graph must correspond to one of a restricted set of identifications of
two tetrahedra along two faces.  Instead of simply selecting rotations
and reflections for each pair of faces one after another and pruning
where possible, we can thus redesign the generation of rotations and
reflections as follows.

\begin{enumerate}
    \item Identify all one-ended chains in the face pairing
    graph.  Recursively select a standard layered solid torus of the
    appropriate size to correspond to each of these one-ended chains.

    \item Identify all remaining double edges in the face pairing graph,
    i.e., all double edges not belonging to the one-ended chains
    previously processed.
    For each such double edge recursively select an identification of
    the two corresponding tetrahedra conforming to
    Lemma~\ref{l-censusdoubleedge}.

    \item At this point many of the individual rotations and reflections
    have already been established.
    Recursively select the remaining rotations and reflections as in
    the original algorithm, pruning where possible.
\end{enumerate}

The first step of this redesigned algorithm should offer a substantial
improvement in running time, as can be seen by the following rough
calculations.  In a one-ended chain on $k$
vertices there are $2k-1$ edges corresponding to $2k-1$ rotations and
reflections
that must be selected.  Using the original algorithm we may be
investigating up to $6^{2k-1} = \frac16 36^k$ possible sets of rotations
and reflections for these edges ($3^{2k-1} = \frac13 9^k$ for orientable
triangulations), although this number will be reduced due to pruning.

Using the redesigned algorithm however, instead of counting all possible
sets of rotations and reflections we need only count the number of
possible standard layered solid tori for this one-ended chain.
In constructing such a layered solid torus we have
two choices for how the base tetrahedron of the layered
solid torus is joined to itself and then two choices for how
each subsequent tetrahedron is layered on (there are in fact three
choices for each layering but one will always lead to an edge of degree
two which from Lemma~\ref{l-pruneedgedeg2} can be ignored).
Thus there are only $2^k$ possible layered solid tori
that can correspond to this one-ended
chain, a vast improvement upon both original estimates of
$\frac16 36^k$ and $\frac13 9^k$ possible sets of rotations and reflections.

\subsection{Statistics} \label{s-statistics}

We close with a summary of statistics illustrating the
effectiveness of the improvements outlined in
Sections~\ref{s-algmeliminatepairings} and~\ref{s-algmredesign}.
Table~\ref{tab-newtimes} lists a series of censuses of closed
minimal {\ppirr} triangulations, where each census is described by the
number of tetrahedra and whether we seek orientable or non-orientable
triangulations.

For each census we present two running times for the required
enumeration of 3-manifold triangulations.  The {\em Old Time} column
contains the running time for the original enumeration algorithm without
the face pairing graph
improvements, and the {\em New Time} column contains the
running time for the algorithm with the improvements outlined above.
Once more all times are displayed as {\em h}:{\em mm}:{\em ss} unless
otherwise indicated, are measured on a single 1.2GHz Pentium~III
processor and were obtained using the program {\regina} \cite{regina}.

\begin{table}[htb]
\begin{center} \begin{tabular}{|c|r|r|r|r|}
    \hline
    & \multicolumn{2}{|c|}{\bf Orientable} &
        \multicolumn{2}{|c|}{\bf Non-Orientable} \\
    \hline
    \bf Tetrahedra &
        \bf Old Time & \bf New Time &
        \bf Old Time & \bf New Time \\
    \hline
    3 & 0:00 & 0:00 & 0:02 & 0:00 \\
    4 & 0:03 & 0:01 & 3:03 & 0:06 \\
    5 & 2:28 & 0:51 & 6:01:53 & 9:19 \\
    6 & 2:49:29 & 56:00 & 5 weeks & 15:08:43 \\ 
    \hline
\end{tabular} \end{center}
\caption{Enumeration statistics with face pairing graph improvements}
\label{tab-newtimes}
\end{table}

Happily we see a marked improvement in running time, particularly for
the enumeration of non-orientable triangulations.  It is indeed due the
algorithm improvements described above that the seven tetrahedron
non-orientable census described in \cite{burton-nor7} was made possible.

\begin{appendix}

\section{Useful Triangulations} \label{a-smalltri}

In Section~\ref{s-prelim} we refer to one-tetrahedron and
two-tetrahedron triangulations of a few specific 3-manifolds.  These
triangulations are presented here for completeness.

\begin{figure}
\centerline{\includegraphics[scale=0.6]{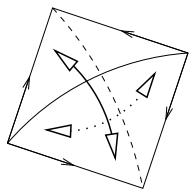}}
\caption{A one-tetrahedron triangulation of $\sss$}
\label{fig-tri3sphere}
\end{figure}

Figure~\ref{fig-tri3sphere} illustrates a one-tetrahedron triangulation
of the 3-sphere $\sss$.  The two front faces are snapped shut about a
degree one edge and the two back faces are similarly snapped shut about
a degree one edge.

\begin{figure}
\centerline{\includegraphics[scale=0.7]{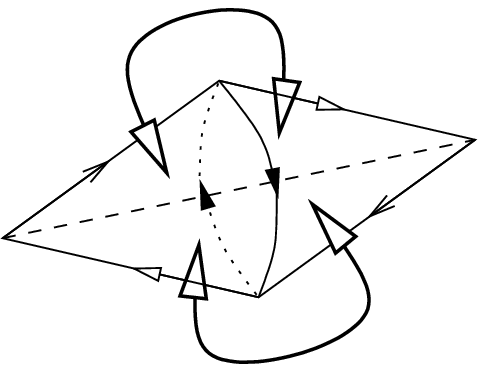}}
\caption{A two-tetrahedron triangulation of $\rps$}
\label{fig-trirps}
\end{figure}

Figure~\ref{fig-trirps} illustrates a two-tetrahedron triangulation of
real projective space $\rps$.  Two tetrahedra are joined along two faces
with an internal edge of degree two, forming a ball with four boundary
faces as illustrated.  The two front faces are then identified with the
two back faces using a $180^\circ$ twist.

\begin{figure}
\centerline{\includegraphics[scale=0.7]{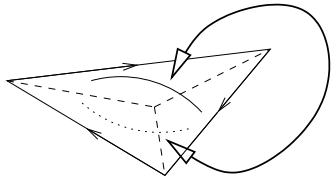}}
\caption{A two-tetrahedron triangulation of $L(3,1)$}
\label{fig-tril31}
\end{figure}

Finally, Figure~\ref{fig-tril31} illustrates a two-tetrahedron
triangulation of the lens space $L(3,1)$.  Two tetrahedra are joined
along three faces with an internal vertex of degree two, forming a
triangular pillow as illustrated.  The front face is then identified
with the back face using a $120^\circ$ twist.

\end{appendix}

\bibliographystyle{amsplain}
\bibliography{topology-20031108}

\vspace{1.5cm}

{\sc \noindent
    \begin{tabbing}
    \hspace{6cm} \= Benjamin A.~Burton \\
    \> Department of Mathematics and Statistics \\
    \> The University of Melbourne, 3010 VIC \\
    \> Australia \\
    \> {\tt bab@debian.org}
    \end{tabbing}
}

\end{document}